\documentclass[12pt,a4paper]{article}
\usepackage[dvipdfmx,hiresbb]{graphicx}
\usepackage{latexsym}
\usepackage{amssymb}
\usepackage{amsmath}
\usepackage{amscd}
\usepackage{mathrsfs}
\usepackage{bm}

\newtheorem{Th}{Theorem}[section]
\newtheorem{Co}[Th]{Corollary}
\newtheorem{Lem}[Th]{Lemma}

\newtheorem{Rem}[Th]{Remark}

\newtheorem{Pro}[Th]{Proposition}
\newcommand{\demo}{\par\noindent{\it Proof. \/}\ }
\newcommand{\enD}{\hfill $\Box$ \vspace{3truemm}\par}
\newcommand{\bx}{\mbox{\boldmath $x$}}
\newcommand{\bX}{\mbox{\boldmath $X$}}

\newcommand{\bv}{\mbox{\boldmath $v$}}
\newcommand{\by}{\mbox{\boldmath $y$}}
\newcommand{\bo}{\mbox{\boldmath $0$}}

\newcommand{\bxi}{\mbox{\boldmath $\xi$}}

\newcommand{\bn}{\mbox{\boldmath $n$}}

\newcommand{\R}{{\mathbb R}}
\newcommand{\lon}{\longrightarrow}

\begin{document}

\title{The theory of graph-like Legendrian unfoldings \\ and its applications}

\author{Shyuichi IZUMIYA}

\date{\today}

\maketitle
\begin{center}
To the memory of my friend Vladimir M. Zakalyukin.
\end{center}
\begin{abstract}
This is mainly a survey article on the recent development of the theory of graph-like Legendrian unfoldings 
and
its applications. The notion of big Legendrian submanifolds was introduced by Zakalyukin for describing the wave front propagations.
Graph-like Legendrian unfoldings belong to a special class of big Legendrian submanifolds.
Although this is a survey article, some new original results and the corrected proofs of some results are given.
\end{abstract}
\renewcommand{\thefootnote}{\fnsymbol{footnote}}
\section{Introduction}
The notion of graph-like Legendrian unfoldings was introduced in \cite{Izumiya93}. 
It belongs to a special class of the big Legendrian submanifolds which Zakalyukin introduced in \cite{Zak76,Zak1}.
There have been some developments on this theory during past two decades\cite{Izumiya93,Izu95,Izumiya-Takahashi,Izumiya-Takahashi2,Izumiya-Takahashi3}.
Most of the results here are already present, implicitly or explicitly, in those articles. 
However, we give in this survey detailed proofs as an aid to understanding and applying the theory. Moreover some of the results here are original, especially Theorem 4.14 which explains 
how the theory of graph-like Legendrian unfoldings is useful for applying to many situations related to the theory of Lagrangian singularities (caustics).
Moreover, it has been known that caustics equivalence (i.e., diffeomorphic caustics) does not imply
Lagrangian equivalence. This is one of the main differences from the theory of Legendrian singularities.
In the theory of Legendrian singularities, wave fronts equivalence (i.e., diffeomorphic wave fronts)
implies Legendrian equivalence generically.
\par
One of the typical examples of big wave fronts (also, graph-like wave fronts) is given by the parallels of a plane curve.
For a curve in the Euclidean plane, its parallels consist of those curves a fixed distance $r$ down the normals in a fixed direction. 
They usually have singularities for sufficiently large $r.$
Their singularities are always Legendrian singularities. 
It is well-known that the singularities of the parallels lie on the evolute of the curve.  
We draw the picture of the parallels of an ellipse and the locus of those singularities in Fig.1.
Moreover, there is another interpretation of the evolute of a curve.
If we consider the family of normal lines to the curve, the evolute is the envelope of this family
of normal lines.
We also draw the envelope of the family of normal lines to an ellipse in Fig.2.
The picture of the corresponding big wave front is depicted in Fig.3.
The evolute is one of the examples of caustics and the family of parallels is a wave front propagation. 
\begin{figure}[htbp]
  \begin{center}
\begin{tabular}{cc}
  \includegraphics[height=35mm]{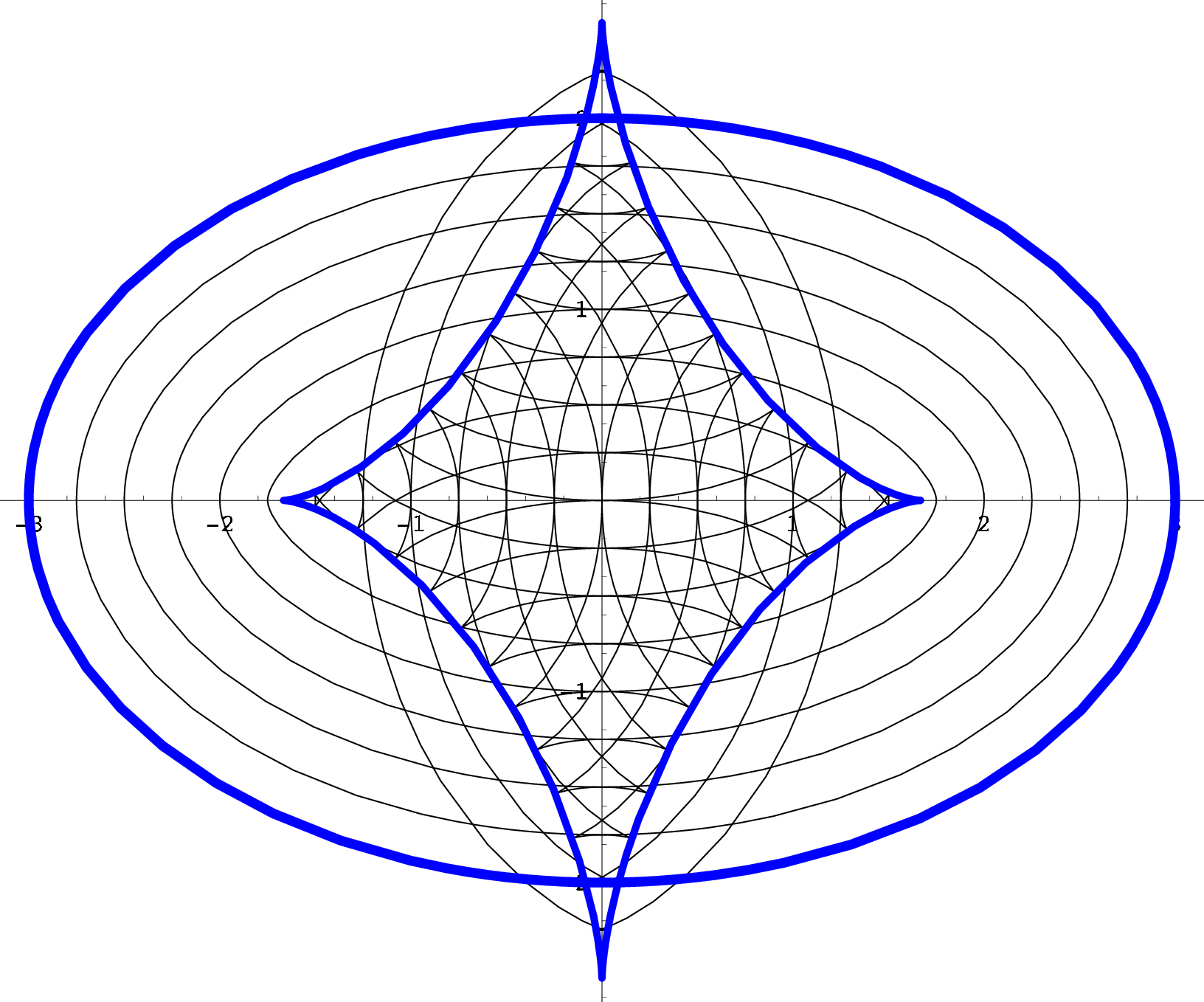} & 
   \includegraphics[height=45mm]{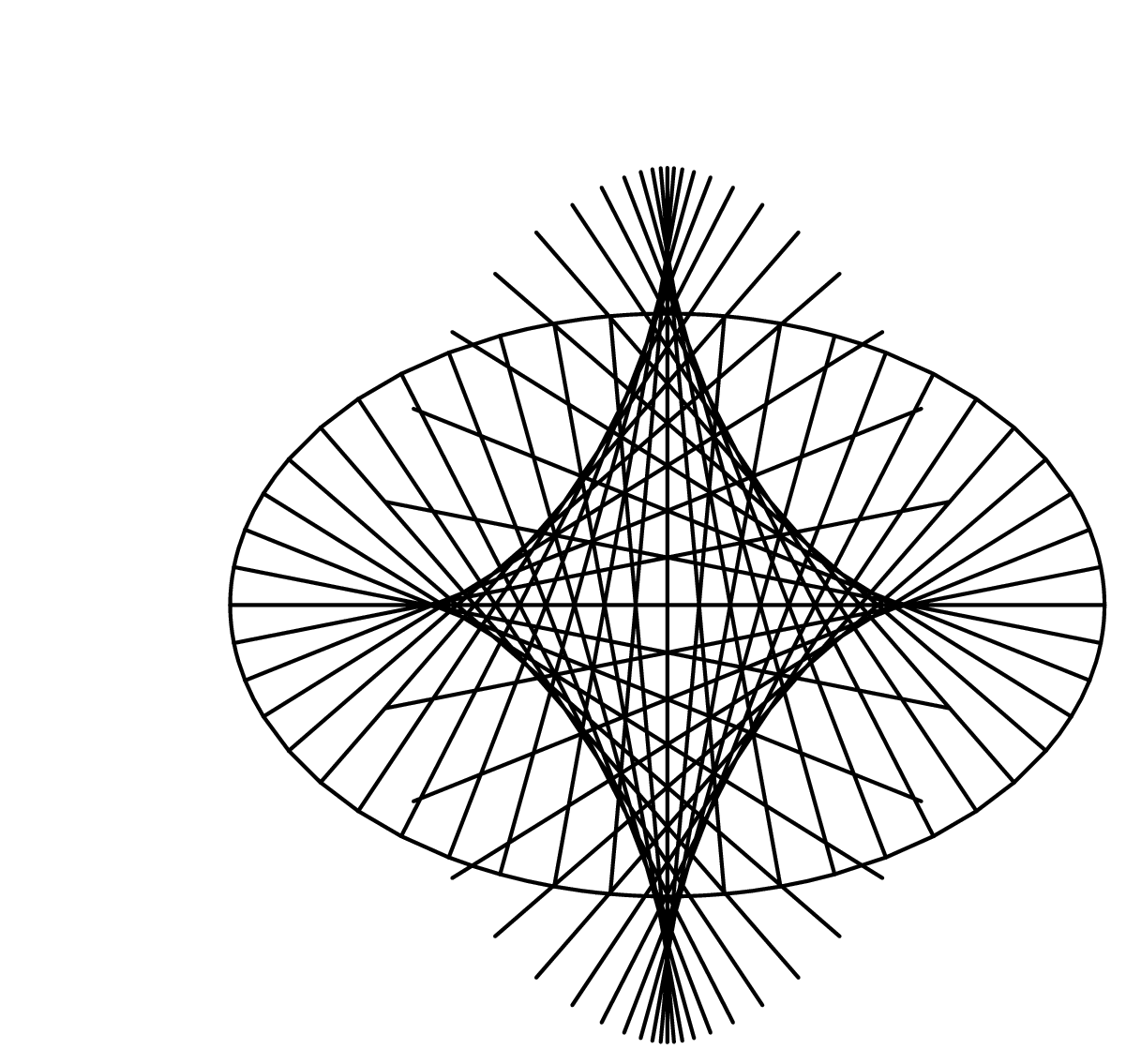}
\\
{Fig.1: The parallels and the evolute} & {Fig.2: The normal lines and the evolute} \\
{of an ellipse} &  {of an ellipse}
\\
\end{tabular}
  \end{center}
\end{figure}

\begin{figure}[htbp]
  \begin{center}
\begin{tabular}{c}
\includegraphics[height=40mm]{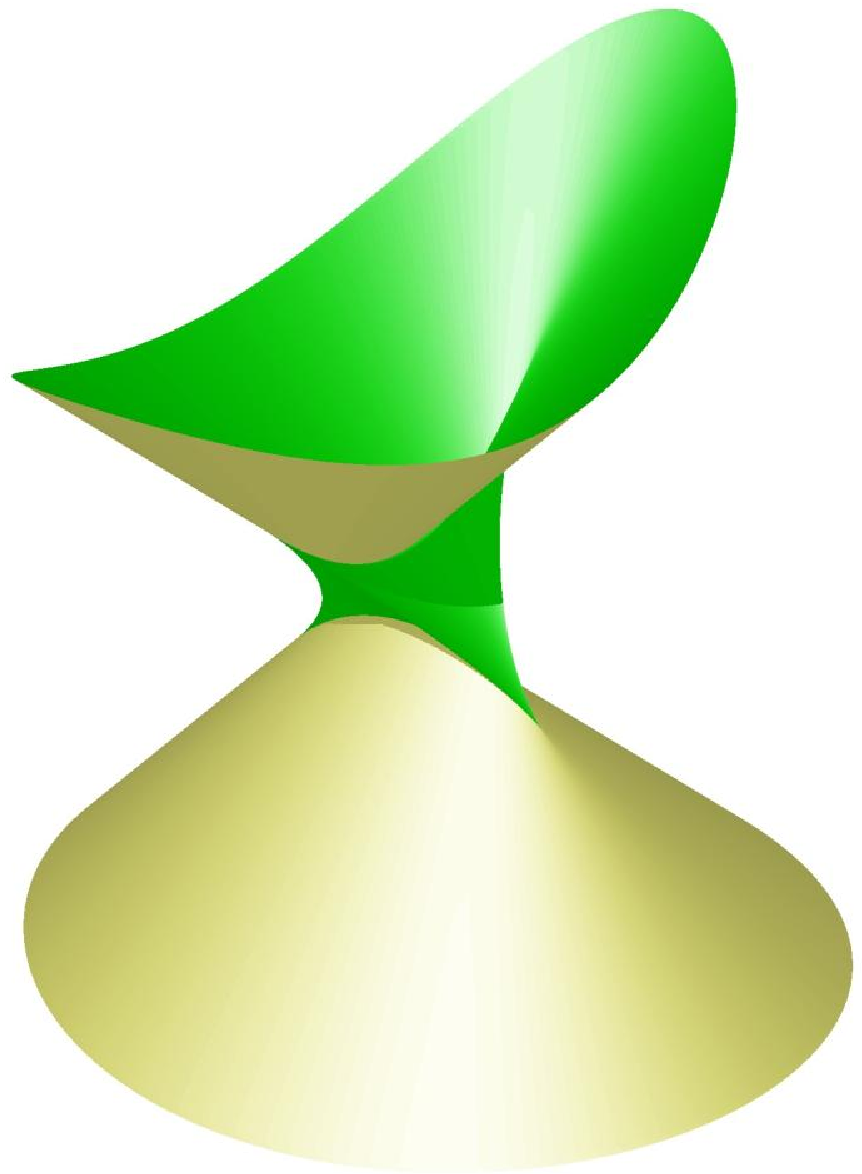} \\
 { Fig.3: The big wave front of the parallels of an ellipse} 
\end{tabular}
  \end{center}
\vskip-.5cm
\end{figure}
\noindent
The caustic is described as the set of critical values of the projection of a Lagrangian submanifold from the phase space onto the configuration space. 
\begin{figure}[htbp]
\vskip-1.5cm
  \begin{center}
\begin{tabular}{c}
\includegraphics[bb=12 132 571 776,clip,height=50mm]{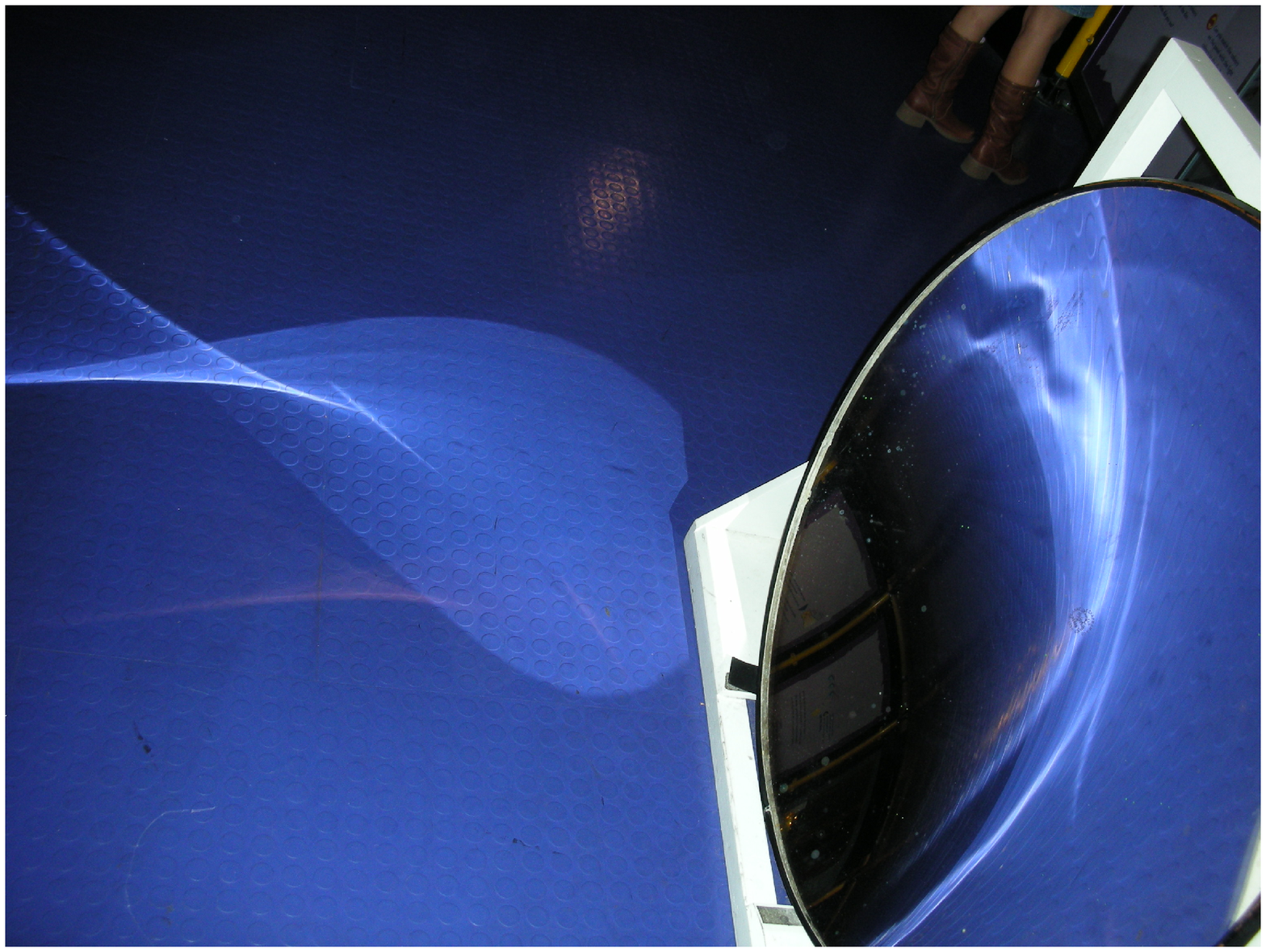} \\
 { Fig.4: The caustic reflected by a mirror} 
\end{tabular}
  \end{center}
\end{figure}
In the real world, the caustics given by reflected rays are visible.
However, the wave front propagations are not visible (cf. Fig. 4). 
Therefore, we can say that there are hidden structures (i.e., wave front 
propagations) 
on the picture of caustics.
In fact, caustics are a subject of classical physics (i.e., optics).
However, the corresponding Lagrangian submanifold is deeply related to the {\it semi-classical approximation} of quantum mechanics (cf., \cite{Hormander,Maslov}).
\par
On the other hand, it was believed around 1989 that the correct framework to describe the parallels of a curve 
is the theory of big wave fronts \cite{Arnold-blue}.
But it was pointed out that $A_1$ and $A_2$ bifurcations do not occur as the parallels of curves \cite{Arnold-pink,Bruce}.
Therefore, the framework of the theory of big wave fronts is too wide for describing the parallels of curves.
The theory of the graph-like Legendrian unfoldings was introduced to construct the correct framework
for the parallels of a curve in \cite{Izumiya93}.
One of the main results in the theory of graph-like Legendrian unfoldings is Theorem 4.14 which
reveals the relation between caustics and wave front propagations.
We give some examples of applications  of the theory of wave front propagations in \S 5.

\section{Lagrangian singularities} 

We give a brief review of the local theory of Lagrangian singularities due to \cite{Arnold1}.
We consider the cotangent bundle $\pi :T^*\R^n\to \R^n$ over $\R^n.$
Let $(x,p)=(x_1,\dots ,x_n,p_1,\dots ,p_n)$ be the canonical coordinates on $T^*\R^n.$
Then the canonical symplectic structure on $T^*\R^n$ is given by the {\it canonical two form}
$\omega =\sum _{i=1}^n dp_i\wedge dx_i.$ 
Let $i:L\subset T^*\R^n$ be a submanifold. We say that $i$ is a {\it Lagrangian submanifold} if
${\rm dim}\, L=n$ and $i^*\omega =0.$ 
In this case, the set of critical values of $\pi\circ i$ is called the {\it caustic} of $i:L\subset T^*\R^n,$
which is denoted by $C_L.$ 
We can interpret the evolute of a plane curve as the caustic of a certain Lagrangian submanifold (cf., \S 5).
One of the main results in the theory of Lagrangian singularities is the description of Lagrangian submanifold germs
by using families of function germs. 
Let $F:(\R^k \times \R^n,0) \to (\R,0)$ 
be an $n$-parameter unfolding of a function germ $f=F|_{\R^k\times \{0\}}:(\R^k,0)\lon (\R,0).$ 
We say that $F$ is a {\it Morse family of functions} if the map germ
$$
\Delta F = \left(\frac{\partial{F}}{\partial{q_1}},\dots,
\frac{\partial{F}}{\partial{q_k}} \right): (\R^k \times \R^n,0) \to (\R^k,0)
$$
is non-singular, 
where $(q,x)=(q_1,\dots,q_k,x_1,\dots,x_n) \in (\R^k \times \R^n,0)$.
In this case, we have a smooth $n$-dimensional submanifold germ 
$C(F)=(\Delta F)^{-1}(0) \subset (\R^k \times \R^n,0)$ and a map germ 
$L(F): (C(F),0) \to T^* \R^n$ defined by
$$
L(F)(q,x) = \left( x, \frac{\partial{F}}{\partial{x_1}}(q,x),\dots,
\frac{\partial{F}}{\partial{x_n}}(q,x) \right).
$$
We can show that $L(F)(C(F))$ is a Lagrangian submanifold germ. 
Then it is known  (\cite{Arnold1}, page 300) that
all Lagrangian submanifold germs in $T^* \R^n$ are constructed by the 
above method.
A Morse family of functions $F:(\R^k \times \R^n,0) \to (\R,0)$ is said to be a {\it generating family} of $L(F)(C(F))$. 
\par
We now define a natural equivalence relation among Lagrangian submanifold germs. 
Let $i: (L,p) \subset (T^* \R^n,p)$ and $i': (L',p') \subset (T^* \R^n,p')$ be 
Lagrangian submanifold germs. Then we say that $i$ and $i'$ are 
{\it Lagrangian equivalent} if there exist a diffeomorphism germ 
$\sigma:(L,p) \to (L',p')$, a symplectic diffeomorphism germ 
$\hat{\tau}:(T^* \R^n,p) \to (T^* \R^n,p')$ and a diffeomorphism germ
$\tau:(\R^n,\pi(p)) \to (\R^n,\pi(p'))$ 
such that $\hat{\tau} \circ i = i' \circ \sigma$ and 
$\pi \circ \hat{\tau} = \tau \circ \pi$, where 
$\pi :(T^* \R^n,p) \to (\R^n,\pi(p))$ is the canonical projection.
Here $\hat{\tau}$ is said to be a {\it symplectic diffeomorphism germ} if it is a diffeomorphism germ such that $\hat{\tau}^*\omega=\omega.$
Then the caustic $C_L$ is diffeomorphic to the caustic $C_{L'}$ by the  diffeomorphism germ $\tau.$
\par
We can interpret the Lagrangian equivalence by using the notion of 
generating families. 
Let $F,G:(\R^k \times \R^n,0) \to (\R,0)$ be function germs. 
We say that $F$ and $G$ are $P$-$\mathcal{R}^+$-{\it equivalent} if 
there exist a diffeomorphism germ 
$\Phi:(\R^k \times \R^n,0) \to (\R^k \times \R^n,0)$ of the form 
$\Phi(q,x)=(\phi_1 (q,x),\phi_2(x))$ and a function germ $h:(\R^n,0) \to (\R,0)$
such that $G(q,x) = F(\Phi(q,x)) + h(x)$. 
For any $F_1:(\R^k\times\R^n,0)\to (\R,0)$ and $F_2 :(\R^{k'}\times\R^n,0)\to (\R,0),$ 
 $F_1$ and $F_2$ are said to be {\it stably} $P$-$\mathcal{R}^+$-{\it equivalent} if they become $P$-$\mathcal{R}^+$-equivalent after 
the addition to the arguments $q_i$ of new arguments $q'_i$ and to the 
functions $F_i$ of non-degenerate quadratic forms $Q_i$ in the new arguments, 
i.e., $F_1 + Q_1$ and $F_2 + Q_2$ are $P$-$\mathcal{R}^+$-equivalent.
Then we have the following theorem:
\begin{Th}\label{Lag.eq}
Let $F:(\R^k\times\R^n,0)\to (\R,0)$ and $G :(\R^{k'}\times\R^n,0)\to (\R,0)$  be Morse families of functions. Then 
$L(F)(C(F))$ and $L(G)(C(G))$ are Lagrangian equivalent if and only if 
$F$ and $G$ are stably $P$-$\mathcal{R}^+$-equivalent.
\end{Th}
\par
Let $F:(\R^k \times \R^n,0) \to (\R,0)$ be a Morse family of functions and
${\cal E}_k$ the ring of function germs of $q=(q_1,\dots ,q_k)$ variables at the origin. 
We say that $L(F)(C(F))$ is {\it Lagrangian stable} if 
$$
\mathcal{E}_k = J_f + {\left\langle 
\frac{\partial{F}}{\partial{x_1}} | _{\R^k \times \{0\}},\dots,
\frac{\partial{F}}{\partial{x_n}} | _{\R^k \times \{0\}}
\right\rangle}_\R +{\langle 1 \rangle}_\R,
$$
where $f=F|_{\R^k\times\{0\}}$ and
$$
 J_f = {\left\langle 
\frac{\partial{f}}{\partial{q_1}}(q),\dots,\frac{\partial{f}}{\partial{q_k}}(q)
\right\rangle}_{\mathcal{E}_k}.
$$
\begin{Rem}{\rm 
In the theory of unfoldings\cite{Bro}, $F$ is said to be an
{\it infinitesimally $P$-$\mathcal{R}^+$-versal unfolding of 
$f=F|_{\R^k \times \{0\}}$} if the above condition is satisfied.
There is a definition of Lagrangian stability (cf., \cite[\S 21.1]{Arnold1}).
It is known that $L(F)(C(F))$ is Lagrangian stable if and only if $F$ is an infinitesimally $P$-$\mathcal{R}^+$-versal unfolding of 
$f=F|_{\R^k \times \{0\}}$ \cite{Arnold1}. In this paper we do not need the original definition of the Lagrangian stability, so that
we adopt the above definition.
}
\end{Rem}

\section{Theory of the wave front propagations}
\par
In this section we give a brief survey of the theory of wave front propagations
(for details, see \cite{Arnold1,Izumiya93,Zak1,Zakalyukin95}, etc).
We consider one parameter families of wave fronts and their bifurcations.
The principal idea is that a one parameter family of wave fronts is considered to be a wave front
whose dimension is one dimension higher than each member of the family.
This is called a {\it big wave front}.
Since the big wave front is a wave front, we start to consider the general theory of Legendrian singularities.
Let $\overline{\pi} :PT^*(\R^m)\lon \R^m$ be the projective cotangent bundle over $\R^m.$
This fibration can be considered as a Legendrian fibration with the canonical contact structure $K$ on $PT^*(\R^m).$
We now review geometric properties of this space.
Consider the tangent bundle
$
\tau :TPT^*(\R^m)\rightarrow PT^*(\R^m)
$
and the differential map
$
d\overline{\pi} :TPT^*(\R^m)\rightarrow T\R^m
$
of $\overline{\pi} .$
For any $X\in TPT^*(\R^m),$ there exists an element
$\alpha\in T^*(\R^m)$ such that
$\tau (X)=[\alpha ].$  For an element $V\in T_x(\R^m),$
the property $\alpha (V)=\bo$ does not depend on the choice of
representative of the class $[\alpha ].$  Thus we can define the canonical
contact structure on $PT^*(\R^m)$ by
$
K=\{X\in TPT^*(\R^m)|\tau (X)(d\overline{\pi} (X))=0\}.
$
We have the trivialization
$
PT^*(\R^m)\cong
\R^m\times P({{\mathbb R}^m}^*)
$
and we call $(x,[\xi])$
{\it homogeneous coordinates}, where
$x=(x_1,\dots ,x_m)\in \R^m$ and
$
[\xi]=[\xi _1:\dots :\xi _m]
$
are homogeneous coordinates of the dual projective space
$P({{\mathbb R}^m}^*).$
It is easy to show that $X\in K_{(x,[\xi])}$ if and only if
$
\sum_{i=1}^m \mu _i\xi _i=0,
$
where
$
d\overline{\pi} (X)=\sum_{i=1}^n \mu _i\frac{\partial}{\partial x_i}.
$
Let $\Phi :(\R^m,0)\lon (\R^m,0)$ be a
diffeomorphism germ. Then we have a unique contact diffeomorphism germ $\widehat{\Phi}:PT^*\R^m\lon PT^*\R^m$
defined by $\widehat{\Phi}(x,[\xi])=(\Phi (x),[\xi\circ d_{\Phi(x)}(\Phi ^{-1})]).$
We call $\widehat{\Phi}$ the {\it contact lift} of $\Phi.$
\par
A submanifold  $i:L\subset PT^*(\R^m)$ is said to be 
a {\it Legendrian submanifold} if $\text{dim}\, L=m-1$ and $di_p(T_pL)\subset K_{i(p)}$
for any $p\in L.$
We also call $\overline{\pi}\circ i=\overline{\pi} |_L: L\lon \R^m$ a {\it Legendrian map} and $W(L)=\overline{\pi} (L)$ a {\it wave front} of $i:L\subset PT^*(\R^m).$
We say that a point $p\in L$ is a {\it Legendrian singular point} if ${\rm rank}\, d(\overline{\pi} |_L)_p < m-1.$
In this case $\overline{\pi} (p)$ is the singular point of $W(L).$
\par
The main tool of the theory of Legendrian singularities is the notion of generating families.
Let $F:({\mathbb R}^k\times{\mathbb R}^m,0 )\longrightarrow ({\mathbb R},0 )$ be a function germ.
We say that $F$ is a {\it Morse family of hypersurfaces} if the map germ
$$
\Delta^*F=\left(F,\frac{\partial F}{\partial q_1},\dots ,\frac{\partial F}{\partial q_k}\right):({\mathbb R}^k\times {\mathbb R}^m,0 )\lon ({\mathbb R}\times {\mathbb R}^k,0 )
$$
is non-singular, where $(q,x)=(q_1,\dots ,q_k,x_1,\dots ,x_m)\in  ({\mathbb R}^k\times{\mathbb R}^m,0 ).$
In this case we have a smooth $(m-1)$-dimensional submanifold germ 
$$
\Sigma _*(F)=\left\{ (q,x)\in ({\mathbb R}^k\times{\mathbb R}^n,0 )\  |\ F(q,x)=\frac{\partial F}{\partial q_1}(q,x)=\cdots =\frac{\partial F}{\partial q_k}(q,x)=0\ \right\}
$$
and we have a map germ $\mathscr{L} _F:(\Sigma _*(F), 0)\lon PT^*{\mathbb R}^m$ defined by 
$$
\mathscr{L} _F(q,x)=\left(x,\left[\frac{\partial F}{\partial x_1}(q,x):\dots :\frac{\partial F}{\partial x_m}(q,x)\right]\right).
$$
We can show that $\mathscr{L}_F(\Sigma _*(F))\subset PT^*(\R^m)$ is a Legendrian submanifold germ. Then it is known (\cite[page 320]{Arnold1}) that all Legendrian submanifold germs in $PT^*({\mathbb R}^m)$ are constructed by the above method.
We call $F$ a {\it generating family} of $\mathscr{L} _F(\Sigma _*(F)).$
Therefore the wave front is given by
$$
W(\mathscr{L} _F(\Sigma _*(F))\! =\!\left\{x\in {\mathbb R}^m\ |\exists q\in {\mathbb R}^k\ {\rm s.t}\ F(q,x)=\frac{\partial F}{\partial q_1}(q,x)=\cdots 
=\frac{\partial F}{\partial q_k}(q,x)=0\ \right\}.
$$
\par
Since the Legendrian submanifold germ $i: (L,p) \subset (PT^*{\mathbb R}^n,p)$ is uniquely determined on the regular part of the wave front $W(L),$
we have the following simple but significant property of Legendrian immersion germs \cite{Zak1}.
\begin{Pro}[Zakalyukin]
Let  $i: (L,p) \subset (PT^*{\mathbb R}^m,p)$  and  
 $i': (L',p') \subset (PT^*{\mathbb R}^m, p')$  be Legendrian immersion germs such that regular sets of $\overline{\pi}\circ i, \overline{\pi}\circ i'$
 are dense respectively.
 Then $(L,p)=(L',p')$ if and only if $(W(L),\overline{\pi}(p))=(W(L'),\overline{\pi} (p'))$.
\end{Pro}
\par
In order to understand the ambiguity of generating families for a fixed Legendrian submanifold germ
we introduce the following equivalence relation among Morse families of hypersurfaces.
Let ${\cal E}_k$ be the local ring of function germs $({\mathbb R}^k,0 )\lon {\mathbb R}$ with 
the unique maximal ideal $\mathfrak{M}_k=\{h\in {\cal E}_k\ |\ h(0)=0\ \}.$
For function germs $F,G : ({\mathbb R}^k\times {\mathbb R}^m,0) 
\lon  ({\mathbb R},0)$,  we say that  $F$  and  $G$  are 
{\it strictly parametrized $\mathcal{K}$-equivalent} (briefly, 
$ S.P$-${\cal K}$-{\it equivalent}) if there exists a diffeomorphism germ  
$\Psi  : ({\mathbb R}^k\times {\mathbb R}^m,0 ) \longrightarrow  
({\mathbb R}^k\times {\mathbb R}^m,0)$ 
of the form  $\Psi (q,x) = (\psi _1(q,x), x)$  for  
$(q,x) \in  ({\mathbb R}^k\times {\mathbb R}^m,0 )$  such that 
$\Psi ^{*}(\langle F\rangle_{ {\cal E}_{k+m}}) = \langle G\rangle_{{ \cal E}_{k+m}}$. Here  $\Psi ^{*} : {\cal E}_{k+m} \lon
{\cal E}_{k+m}$  is the pull back ${\mathbb R}$-algebra isomorphism defined by  
$\Psi ^{*}(h) =  h\circ \Psi.$
The definition of {\it stably $S.P$-$\mathcal{K}$-equivalence} among Morse families of hypersurfaces is similar to
the definition of stably $P$-$\mathcal{R}^+$-equivalence among Morse families of functions.
The following is the key lemma of the theory of Legendrian singularities (cf. \cite{Arnold1,Gor-Zak,Zak}).
\begin{Lem}[Zakalyukin]
Let $F:(\R^k\times\R^m,0)\to (\R,0)$ and $G :(\R^{k'}\times\R^m,0)\to (\R,0)$  be Morse families of hypersurfaces. Then 
$(\mathscr{L}_F(\Sigma _*(F)),p)=(\mathscr{L}_G(\Sigma _*(G)),p)$ if and only if 
$F$ and $G$ are stably $S.P$-$\mathcal{K}$-equivalent.
\end{Lem}
\par
Let $F:(\R^k\times\R^m,0)\to (\R,0)$ be a Morse family of hypersurfaces and $\Phi:(\R^m,0)\lon (\R^m,0)$ a
diffeomorphism germ. We define $\Phi^*F:(\R^k\times\R^m,0)\to (\R,0)$ by
$\Phi^*F(q,x)=F(q,\Phi (x)).$ Then we have $(1_{\R^q}\times \Phi)(\Sigma _*(\Phi^*F))=\Sigma _*(F)$
and 
\[
\mathscr{L}_{\Phi^*F}(\Sigma _*(\Phi^*F))=\left\{ \left(x, \left[\left(\frac{\partial F}{\partial x}(q,\Phi(x))\right)\circ d\Phi _x\right]\right)\ \Bigm | (q,\Phi(x))\in \Sigma _*(F)\right\},
\]
so that $\widehat{\Phi}(\mathscr{L}_{\Phi^*F}(\Sigma _*(\Phi^*F)))=\mathscr{L}_F(\Sigma _*(F))$ as set germs.
\begin{Pro}
Let $F:(\R^k\times\R^m,0)\to (\R,0)$ and $G :(\R^{k'}\times\R^m,0)\to (\R,0)$  be Morse families of hypersurfaces.
For a diffeomorphism germ $\Phi:(\R^m,0)\lon (\R^m,0),$ $\widehat{\Phi}(\mathscr{L}_G(\Sigma _*(G)))=\mathscr{L}_F(\Sigma _*(F))$ if and only if $\Phi^*F$ and $G$ are stably $S.P$-$\mathcal{K}$-equivalent.
\end{Pro}
\demo
Since $\widehat{\Phi}(\mathscr{L}_{\Phi^*F}(\Sigma _*(\Phi^*F)))=\mathscr{L}_F(\Sigma _*(F)),$ we have
$\mathscr{L}_{\Phi^*F}(\Sigma _*(\Phi^*F))=\mathscr{L}_G(\Sigma _*(G)).$
By Lemma 3.2, the assertion holds.
\enD
We say that $\mathscr{L}_F(\Sigma _*(F))$ and $\mathscr{L}_G(\Sigma _*(G))$ are {\it Legendrian equivalent}
if there exists a diffeomorphism germ $\Phi:(\R^m,0)\lon (\R^m,0)$ such that the condition in the above proposition holds. By Lemma 3.1, under the generic condition on $F$ and $G$, 
 $\Phi(W(\mathscr{L}_G(\Sigma _*(G))))=W(\mathscr{L}_F(\Sigma _*(F)))$ if and only if $\widehat{\Phi}(\mathscr{L}_G(\Sigma _*(G)))=\mathscr{L}_F(\Sigma _*(F))$ for a diffeomorphism germ $\Phi:(\R^m,0)\lon (\R^m,0)$.
\par
We now consider the case when $m=n+1$ and distinguish space and time coordinates, so that 
we denote that $\R^{n+1}=\R^n\times \R$ and coordinates
are denoted by $(x,t)=(x_1,\dots, x_n,t)\in \R^n\times \R.$
Then we consider the
projective cotangent 
bundle $
\overline{\pi} :PT^*(\R^n\times\R)\to\R^n\times\R.
$ 
Because of the trivialization
$
PT^*(\R^n\times\R)\cong
(\R^n\times\R)\times P(({\mathbb R}^{n}\times \R)^*),
$
we have homogeneous coordinates
$
((x_1,\dotsc ,x_n,t),[\xi _1:\cdots:\xi _n:\tau ]).
$
We remark that $PT^*(\R^n\times\R)$ is a fiber-wise compactification of the 
1-jet space as follows: 
We consider an affine open subset $U_\tau =\{((x,t),[\xi :\tau ])|\tau\not= 0\}$ of $PT^*(\R^n\times\R).$
For any $((x,t),[\xi :\tau ])\in U_\tau ,$ we have
$$
((x_1,\dots ,x_n,t),[\xi _1:\cdots:\xi _n:\tau ])
=((x_1,\dots ,x_n,t),[-({\xi _1 }/{\tau }):\cdots :-({\xi _n }/{\tau }):-1 ]),
$$
so that we may adopt the corresponding {\it affine coordinates}
$
((x_1,\dots ,x_n,t),(p_1,\dots ,p_n)),
$
where $p_i=-{\xi _i}/{\tau }.$
On $U_\tau$ we can easily show that $\theta ^{-1}(0)=K|_{U_\tau} ,$
where $\theta =dt-\sum_{i=1}^n p_idx_i.$
This means that $U_\tau $ can be identified with the 1-jet space which is denoted by
$J^1_{GA}(\R^n,\R)\subset PT^*(\R^n\times\R).$
We call the above coordinates {\it a system of graph-like affine coordinates.}
Throughout this paper, we use this identification.
\par
For a Legendrian submanifold $i:L\subset PT^*(\R^n\times\R),$
the corresponding wave front $\overline{\pi}\circ i(L)=W(L)$ is called a {\it big wave front.}
We call 
$$
W_t(L)=\pi_1(\pi _2^{-1}(t) \cap W(L))\quad (t\in \R)
$$
a {\it momentary front} (or, a {\it small front}) for each $t\in (\R,0),$
where $\pi_1:\R^n \times \R \to \R^n$ and $\pi_2:\R^n \times \R \to \R$ 
are the canonical projections defined by $\pi_1(x,t)=x$ and $\pi_2(x,t)=t$ respectively. 
In this sense, we call $L$ a {\it big Legendrian submanifold.}
We say that a point $p\in L$ is a {\it space-singular point} if ${\rm rank}\, d(\pi_1\circ \overline{\pi} |_L)_p < n$
and a {\it time-singular point} if ${\rm rank}\, d(\pi _2\circ\overline{\pi} |_L)_p=0,$ respectively.
By definition, if $p\in L$ is a Legendrian singular point, then it is a space-singular point of $L.$
Even if we have no Legendrian singular points, we have space-singular points. In this case we have the
following lemma.
\begin{Lem} Let $i:L\subset PT^*(\R^n\times\R)$ be a big Legendrian submanifold without Legendrian singular points.
If $p\in L$ is a space-singular point of $L$, then $p$ is not a time-singular point of $L.$
\end{Lem}
\demo
By the assumption, $\overline{\pi}|_L$ is an immersion. For any $v\in T_pL,$ there exists $X_v\in T_{\overline{\pi} (p)}(\R^n\times \{0\})$ and
$Y_v\in T_{\overline{\pi} (p)}(\{0\}\times \R)$ such that $d(\overline{\pi}|_L)_{p}(v)=X_v+Y_v.$
If ${\rm rank}\, d(\pi _2\circ\overline{\pi} |_L)_p=0,$ then $d(\overline{\pi}|_L)_{p}(v)=X_v$ for any $v\in T_pL.$
Since $p$ is a space-singular point of $L$, there exits $v\in T_pL$ such that $X_v=0,$ so that $d(\overline{\pi} |_L)_p(v)=0.$
This contradicts to the fact that $\overline{\pi} |_L$ is an immersion.
\enD
\par
The {\it discriminant of the family $\{W_t(L)\}_{t\in (\R,0)}$} is defined as the image of singular points of $\pi _1|_{W(L)}.$
In the general case, the discriminant consists of three components: {\it the caustic} $C_L=\pi_1(\Sigma (W(L))$, where $\Sigma (W(L))$ is the set of singular points of $W(L)$ (i.e, the critical value set of the Legendrian mapping $\overline{\pi}|_{L}$),
{\it the Maxwell stratified set} $M_L,$ the projection of the closure of the self intersection set of $W(L);$
and also of the critical value set $\Delta_L$ of $\pi_1 |_{W(L)\setminus \Sigma (W(L))}.$ In \cite{Izumiya-Takahashi,Zakalyukin95}, it has been stated that $\Delta_L$ is the {\it envelope of the family of momentary fronts}. However, we remark that $\Delta_L$ is not necessarily the envelope of the family of the projection of smooth momentary fronts $\overline{\pi}(W_t(L)).$
It can happen that $\pi _2^{-1}(t)\cap W(L)$ is non-singular but $\pi _1|_{\pi _2^{-1}(t)\cap W(L)}$ has singularities,
so that $\Delta_L$ is the set of critical values of the family of mappings $\pi _1|_{\pi _2^{-1}(t)\cap W(L)}$ for
smooth $\pi _2^{-1}(t)\cap W(L)$ (cf., \S 5.2). 
\par
For any Legendrian submanifold germ $i:(L,p_0)\subset  (PT^*(\R^n\times \R),p_0),$ 
there exists a generating family. 
Let
$
\mathcal{F}:(\R^k\times(\R^n\times\R),0)\to (\R,0)
$
be a Morse family of hypersurfaces.
In this case, we call $\mathcal{F}$ a {\it big Morse family of hypersurfaces.}  Then $\Sigma_*(\mathcal{F})=\Delta ^*(\mathcal{F})^{-1}(0)$ is
a smooth $n$-dimensional submanifold germ. 
By the previous arguments, we have a big Legendrian submanifold
$\mathscr{L}_{\mathcal{F}}(\Sigma _*(\mathcal{F}))$ where
$$
\mathscr{L} _\mathcal{F}(q,x,t)=\left(x,t,\left[\frac{\partial \mathcal{F}}{\partial x}(q,x,t):\frac{\partial \mathcal{F}}{\partial t}(q,x,t)\right]\right),
$$
and
$$
\left[\frac{\partial \mathcal{F}}{\partial x}(q,x,t):\frac{\partial \mathcal{F}}{\partial t}(q,x,t)\right]
=\left[\frac{\partial \mathcal{F}}{\partial x_1}(q,x,t):\cdots :\frac{\partial \mathcal{F}}{\partial x_n}(q,x,t):
\frac{\partial \mathcal{F}}{\partial t}(q,x,t)\right].
$$
\par
We now consider an equivalence relation among big Legendrian submanifolds which preserves the discriminant of families of momentary fronts.
The following equivalence relation among big Legendrian submanifold germs has been independently introduced in 
\cite{Izu95, Zakalyukin95} for different purposes:
Let $i:(L,p_0)\subset (PT^*(\R^n\times\R),p_0)$ and $i':(L',p_0')\subset (PT^*(\R^n\times\R),p_0')$
be big Legendrian submanifold germs. 
We say that $i$ and $i'$ are {\it strictly parametrized$^+$ Legendrian equivalent} (or, briefly {\it $S.P^+$-Legendrian equivalent}) 
if there exists a diffeomorphism germs $\Phi :(\R^n\times\R ,\overline{\pi}(p_0))\to (\R^n\times\R ,\overline{\pi}(p_0'))$
of the form $\Phi (x,t)=(\phi _1(x),t+\alpha (x))$ such that
$\widehat{\Phi}(L)=L'$ as set germs, where 
$\widehat{\Phi}:(PT^*(\R^n\times\R),p_0) \to (PT^*(\R^n\times\R),p'_0)$ is the unique contact lift of $\Phi$.
We can also define the notion of stability of Legendrian submanifold germs with respect to $S.P^+$-Legendrian equivalence which is analogous to the stability of Lagrangian submanifold germs with respect to Lagrangian equivalence (cf. [1, Part III]).
We investigate $S.P^+$-Legendrian equivalence by using the notion of generating families of Legendrian submanifold germs.
Let $\overline{f},\overline{g}:(\R^k\times\R,0) \to (\R,0)$ be function germs.  Remember that
$\overline{f}$ and $\overline{g}$ are S.P-${\cal K}$-equivalent if
there exists a diffeomorphism germ
$
\Phi :(\R^k \times \R,0)\to (\R^k \times \R,0)
$
of the form $\Phi (q,t)=(\phi (q,t),t)$ such that
$
\langle \overline{f}\circ\Phi\rangle _{{\cal E}_{k+1}}=
\langle \overline{g}\rangle _{{\cal E}_{k+1}}.
$
Let
$\mathcal{F},\mathcal{G}:(\R^k\times (\R^n\times\R),0)\to (\R,0)$ be function
germs.  We say that $\mathcal{F}$ and $\mathcal{G}$ are {\it space-$S.P^+$-${\cal K}$-equivalent} (or, briefly, {\it $s$-$S.P^+$-${\cal K}$-equivalent}) if there exists a diffeomorphism germ
$
\Psi :(\R^k\times (\R^n\times\R),0) \to (\R^k\times (\R^n\times\R),0)
$
of the form 
$
\Psi (q,x,t)=(\phi (q,x,t),\phi_1(x),t+\alpha (x))
$
such that
$
\langle F\circ\Psi \rangle _{{\cal E}_{k+n+1}}=
\langle G\rangle _{{\cal E}_{k+n+1}}.
$
The notion of $S.P^+$-${\cal K}$-versal deformation
plays an important role for our purpose. 
 We define the extended tangent space of $\overline{f}:(\R^k\times \R,0)\to (\R,0)$
 relative to $S.P^+$-${\cal K}$ by
$$
T_e(S.P^+\mbox{\rm -}{\cal K})(\overline{f})=\left\langle \frac{\partial \overline{f}}{\partial q_1},\dots, \frac{\partial \overline{f}}{\partial q_k},\overline{f} \right\rangle _{{\cal E}_{k+1}}+
\left\langle \frac{\partial \overline{f}}{\partial t} \right\rangle _{\R}.
$$
Then we say that $F$ is an {\it infinitesimally $S.P^+$-${\cal K}$-versal} deformation of $\overline{f}=F|_{\R^k \times \{0\} \times \R }$ if it satisfies
$$
{\cal E}_{k+1}=T_e(S.P^+\mbox{\rm -}{\cal K})(\overline{f})+
\left\langle \frac{\partial \mathcal{F}}{\partial x_1}|_{\R^k\times\{0\}\times \R} ,\dots
,\frac{\partial \mathcal{F}}{\partial x_n}|_{\R^k\times\{0\}\times\R} \right\rangle _{\R}.
$$

\begin{Th}{\rm \cite{Izu95, Zakalyukin95}}
Let $\mathcal{F}:(\R^k\times(\R^n\times\R),0)\to (\R,0)$
and $\mathcal{G}:(\R^{k'}\times (\R^n\times\R),0)\to (\R,0)$
be big Morse families of hypersurfaces.  Then
\par\noindent
	{\rm ($1$)} $\mathscr{L} _\mathcal{F}(\Sigma _*(\mathcal{F}))$ and $\mathscr{L}_\mathcal{G}(\Sigma _*(\mathcal{G}))$ are $S.P^+$-Legendrian equivalent if and only if
$\mathcal{F}$ and $\mathcal{G}$ are stably $s$-$S.P^+$-${\cal K}$-equivalent.
\par\noindent
{\rm ($2$)} $\mathscr{L} _\mathcal{F}(\Sigma _*(\mathcal{F}))$ is $S.P^+$-Legendre stable if and only if $\mathcal{F}$ is an infinitesimally $S.P^+$-${\cal K}$-versal deformation of $\overline{f}=\mathcal{F}|_{\R^k\times\{0\}\times \R}.$
\end{Th}
\demo
By definition, $\mathcal{F}$ and $\mathcal{G}$ are stably $s$-$S.P^+$-${\cal K}$-equivalent if there exists
a diffeomorphism germ $
\Phi :(\R^n\times\R,0) \to (\R^n\times\R,0)
$
of the form 
$
\Phi (x,t)=(\phi_1(x),t+\alpha (x))
$
such that $\Phi^*F$ and $G$ are stably $S.P$-$\mathcal{K}$-equivalent.
By Proposition 3.3, we have the assertion (1).
For the proof of the assertion (2), we need some more preparations, so that we omit it.
We only remark here that the proof is analogous to the proof of \cite[Theorem in \S 21.4]{Arnold1}.
\enD
The assumption in Proposition 3.1 is a generic condition for $i,i'.$
Especially, if $i$ and $i'$ are $S.P^+$-Legendre stable, then these big Legendrian submanifold germs satisfy the assumption. 
Concerning the discriminant and the bifurcation of momentary fronts, we define the following equivalence relation among big wave front germs.
Let $i: (L,p_0) \subset (PT^*({\mathbb R}^n\times\R),p_0)$  and  
 $i': (L',p'_0) \subset (PT^*({\mathbb R}^n\times\R), p'_0)$  be big Legendrian submanifold germs. 
We say that $W(L)$ and $W(L')$ are {\it $S.P^+$-diffeomorphic} if there exists 
a diffeomorphism germ $\Phi :(\R^n\times\R ,\overline{\pi} (p_0))\to (\R^n\times\R ,\overline{\pi} (p'_0))$
of the form $\Phi (x,t)=(\phi _1(x),t+\alpha (x))$ such that $\Phi (W(L))=W(L')$ as set germs. 
We remark that $S.P^+$-diffeomorphism among big wave front germs preserves the diffeomorphism types of the discriminants \cite{Zakalyukin95}.
By Proposition 3.1, we have the following proposition.

\begin{Pro}
Let  $i: (L,p_0) \subset (PT^*({\mathbb R}^n\times\R),p_0)$  and  
 $i': (L',p_0') \subset (PT^*({\mathbb R}^n\times\R), p_0')$  be big Legendrian submanifold germs such that regular sets of $\overline{\pi}\circ i, \overline{\pi}\circ i'$
 are dense respectively.
 Then $i$ and $i'$ are $S.P^+$-Legendrian equivalent  if and only if $(W(L),\overline{\pi} (p_0))$ and $( W(L'),\overline{\pi} (p'_0))$ are $S.P^+$-diffeomorphic.
 \end{Pro}

\begin{Rem} {\rm If we consider a diffeomorphism germ $
\Phi :(\R^n\times\R,0) \to (\R^n\times\R,0)
$
of the form 
$
\Phi (x,t)=(\phi_1(x,t), \phi _2 (t)),
$
we can define a {\it time-Legendrian equivalence} among big Legendrian submanifold germs.
We can also define a {\it time-$P$-$\mathcal{K}$-equivalence} among big Morse families of hypersurfaces.
By the similar arguments to the above paragraphs, we can show that these equivalence relations describe
the bifurcations of momentary fronts of big Legendrian submanifolds.
In \cite{Zak1} Zakalyukin classified generic big Legendrian submanifold germs by time-Legendrian equivalence.
The natural equivalence relation among big Legendrian submanifold germs is induced by
diffeomorphism germs $
\Phi :(\R^n\times\R,0) \to (\R^n\times\R,0)
$
of the form 
$
\Phi (x,t)=(\phi_1(x), \phi _2 (t)).
$
This equivalence relation classifies both the discriminants and the bifurcations of momentary fronts of big Legendrian submanifold germs.
However, it induces an equivalence relation among divergent diagrams $(\R^n,0)\leftarrow (\R^{n+1},0)\rightarrow (\R,0),$
so that it is almost impossible to have a classification by this equivalence relation.
Here, we remark that the corresponding group of the diffeomorphisms is not a geometric subgroup of $\mathcal{A}$ and $\mathcal{K}$
in the sense of Damon\cite{Damon}.
Moreover, if we consider a diffeomorphism germ $
\Phi :(\R^n\times\R,0) \to (\R^n\times\R,0)
$
of the form 
$
\Phi (x,t)=(\phi_1(x), t),
$
we have a stronger equivalence relation among big Legendrian submanifolds,
which is called an {\it $S.P$-Legendrian equivalence}. Although this equivalence relation gets rid of the difficulty for the above equivalence relation, there appear function moduli for generic classifications in very low dimensions (cf., \S 5).
In order to avoid the function moduli, we introduced the $S.P^+$-Legendrian equivalence.
If we have a generic classification of big Legendrian submanifold germs by $S.P^+$-Legendrian equivalence, we have a classification by the $S.P$-Legendrian equivalence modulo function moduli.
See \cite{Izu95, Zakalyukin95} for details.
\par
On the other hand, we can also define a {\it space-Legendrian equivalence} among big Legendrian submanifold germs.
According to the above paragraphs, we use a diffeomorphism germ $\Phi:(\R^n\times\R,0)\lon (\R^n\times\R,0)$
of the form $\Phi (x,t)=(\phi _1(x),\phi_2 (x,t)).$
The corresponding equivalence among big Morse families of hypersurfaces is the {\it space-$P$-$\mathcal{K}$-equivalence}
which is analogous to the above definitions (cf., \cite{Izudoc}). Recently, we discovered an application of this equivalence relation
to the geometry of world sheets in Lorentz-Minkowski space. See \cite{Izu14} for details.
}
\end{Rem}

\section{Graph-like Legendrian unfoldings}
In this section we explain the theory of 
graph-like Legendrian unfoldings. 
Graph-like Legendrian unfoldings belong to a special class of big Legendrian 
submanifolds. 
A big Legendrian submanifold $i:{L}\subset PT^*(\R^n\times\R)$ is said to be a {\it graph-like Legendrian unfolding} if $L\subset J^1_{GA}(\R^n,\R)$. 
We call $W(L)=\overline{\pi}(L)$ a {\it graph-like wave front} of $L,$ where $\overline{\pi}:J^1_{GA}(\R^n,\R)\lon \R^n\times\R$
is the canonical projection.
We define a mapping $\Pi:J^1_{GA}(\R^n,\R)\lon T^*\R^n$ by
$\Pi(x,t,p)=(x,p),$ where $(x,t,p)=(x_1,\dots, x_n,t,p_1,\dots ,p_n)$ and the canonical contact form
on $J^1_{GA}(\R^n,\R)$ is given by $\theta =dt-\sum_{i=1}^n p_idx_i.$
Here, $T^*\R^n$ is a symplectic manifold with the
canonical symplectic structure $\omega=\sum_{i=1}^n dp_i\wedge dx_i$ (cf. \cite{Arnold1}).
Then we have the following proposition.
\begin{Pro}[\cite{Izumiya-Takahashi}]
For a graph-like Legendrian unfolding $L\subset J^1_{GA}(\R^n,\R),$ $z\in L$ is a singular point of $\overline{\pi}|_L: L\lon \R^n\times\R$
if and only if it is a singular point of $\pi_1\circ\overline{\pi}|_L:L\lon \R^n.$
Moreover, $\Pi|_L:L\lon T^*\R^n$ is immersive, so that $\Pi(L)$ is a Lagrangian submanifold in $T^*\R^n.$
\end{Pro}
\demo
Let $z\in L$ be a singular point of $\pi_1\circ\overline{\pi}|_L.$ Then there exists a non-zero tangent vector $\bv\in T_zL$ such that $d(\pi_1\circ\overline{\pi}|_L)_z(\bv)=0.$ For the canonical coordinate $(x,t,p)$ of $J^1_{GA}(\R^n,\R),$
we have
\[
\bv=\sum_{i=1}^n\alpha _i\frac{\partial}{\partial x_i}+\beta\frac{\partial}{\partial t}+\sum_{j=1}^n\gamma _j\frac{\partial}{\partial p_j}
\]
for some real numbers $\alpha _i,\beta,\gamma _j\in \R.$
By the assumption, we have $\alpha _i=0$ $(i=1,\dots ,n).$
Since $L$ is a Legendrian submanifold in $J^1_{GA}(\R^n,\R),$
we have $0=\theta(\bv)=\beta -\sum_{i=1}^n\gamma _i\alpha _i=\beta.$
Therefore, we have
\[
d\overline{\pi}(\bv)=\sum_{i=1}^n\alpha _i\frac{\partial}{\partial x_i}+\beta\frac{\partial}{\partial t}=\bo.
\]
This means that $z\in L$ is a singular point of $\overline{\pi}|_L.$
The converse assertion holds by definition.
\par
We consider a vector $\bv\in T_zL$ such that $d\Pi_z(\bv)=\bo.$
For similar reasons to the above case, we have $\bv=\bo.$
This means that $\Pi|_L$ is immersive.
Since $L$ is a Legendrian submanifold in $J^1_{GA}(\R^n,\R),$ we have
\[
\omega|_{\Pi(L)}=(\Pi|_L)^*\omega=\Pi^*\omega|_L=d\theta |_L=d(\theta |_L)=0.
\]
This completes the proof.
\enD
We have the following corollary of Proposition 4.1.
\begin{Co}
For a graph-like Legendrian unfolding $L\subset J^1_{GA}(\R^n,\R),$ $\Delta _L$ is
the empty set, so that the discriminant of the family of momentary fronts is $C_L\cup M_L.$
\end{Co}
Since ${L}$ is a big Legendrian submanifold in $PT^*(\R^n\times\R),$
it has a generating family ${\mathcal F} :(\R^k\times (\R^n\times\R),0)\to (\R,0)$ at least locally.
Since $L\subset J^1_{GA}(\R^n,\R)=U_\tau\subset PT^*(\R^n\times \R),$
it satisfies the condition $(\partial {\mathcal F}/\partial t)(0)\not= 0.$
Let ${\mathcal F} :(\R^k\times (\R^n\times\R),0)\to (\R,0)$ be a big Morse family of hypersurfaces.
We say that ${\mathcal F}$ is a {\it graph-like Morse family of hypersurfaces}
if $(\partial {\mathcal F}/\partial t)(0)\not= 0.$
It is easy to show that the corresponding big Legendrian submanifold germ
is a graph-like Legendrian unfolding. 
Of course, all graph-like Legendrian unfolding germs can be constructed by the above way.
We say that ${\mathcal F}$ is a {\it graph-like generating family} of
$\mathscr{L}_{\mathcal F}(\Sigma _*({\mathcal F})).$
We remark that the notion of graph-like Legendrian unfoldings and corresponding generating families have been introduced in \cite{Izumiya93} to describe the perestroikas of wave fronts
given as the solutions for general eikonal equations. In this case, there is an additional condition. We say that 
$\mathcal F:(\R^k\times (\R^n\times\R),0)\to (\R,0)$ is {\it non-degenerate} if $\mathcal F$ satisfies the conditions $(\partial {\mathcal F}/\partial t)(0)\not= 0$ and 
$\Delta^*\mathcal{F}|_{\R^k \times \R^n \times \{0\}}$ is a submersion germ. In this case we call $\mathcal F$ a {\it non-degenerate graph-like generating family}.
We have the following proposition.
\begin{Pro} Let $\mathcal{F}:(\R^k\times (\R^n\times\R),0)\to (\R,0)$ be a graph-like Morse family of hypersurfaces. Then $\mathcal{F}$ is non-degenerate if and only if $\pi_2\circ\overline{\pi}|_{\mathscr{L}_{\mathcal{F}}(\Sigma_{*}(\mathcal{F}))}$ is submersive.
\end{Pro}
\demo
By the definition of $\mathscr{L}_{\mathcal{F}},$ we have $\pi_2\circ\overline{\pi}|_{\mathscr{L}_{\mathcal{F}}(\Sigma_{*}(\mathcal{F}))}=\pi _2\circ\pi _{n+1}|_{\Sigma _*(\mathcal{F})},$ where $\pi _{n+1}:\R^k\times \R^n\times \R\lon \R^n\times \R$
is the canonical projection.
Since $\Sigma _*(\mathcal{F})=\Delta ^*(\mathcal{F})^{-1}(0)\subset (\R^k\times (\R^n\times \R),0),$
$\pi _2\circ\pi _{n+1}|_{\Sigma _*(\mathcal{F})}$ is submersive if and only if
\[
{\rm rank}\,\left(\frac{\partial \Delta ^*(\mathcal{F})}{\partial q}(0),\frac{\partial \Delta ^*(\mathcal{F})}{\partial x}(0)\right)=k+1.
\]
The last condition is equivalent to the condition that
\[
\Delta ^*(\mathcal{F}|_{\R^k\times\R^n\times \{0\}}):(\R^k\times \R^n\times \{0\},0)\lon (\R\times \R^k,0)
\]
is non-singular.
This completes the proof.
\enD
We say that a graph-like Legendrian unfolding $L\subset J^1_{GA}(\R^n,\R)$ is {\it non-degenerate} if
$\pi _2\circ\overline{\pi}|L$ is submersive.
The notion of graph-like Legendrian unfolding was introduced in \cite{Izumiya93}.
Non-degeneracy was then assumed for general graph-like Legendrian unfoldings.
However, during the last two decades, we have clarified the situation and
non-degeneracy is now defined as above.
\par
We can consider the following more restrictive class of graph-like generating families: 
Let ${\mathcal F}$ be a graph-like Morse family of hypersurfaces.
By the implicit function theorem, there exists a function
$F:(\R^k\times \R^n,0)\to (\R,0)$
such that $\langle {\mathcal F}(q,x,t)\rangle_{{\mathcal E}_{k+n+1}}
=\langle F(q,x)-t\rangle_{{\mathcal E}_{k+n+1}}.$
Then we have the following proposition.
\begin{Pro}
Let $\mathcal{F}:(\R^k\times (\R^n\times\R),0)\to (\R,0)$ and $F:(\R^k\times \R^n,0)\to (\R,0)$ be function germs such that 
$\langle {\mathcal F}(q,x,t)\rangle_{{\mathcal E}_{k+n+1}}
=\langle F(q,x)-t\rangle_{{\mathcal E}_{k+n+1}}.$ Then $\mathcal{F}$ is a graph-like Morse family of hypersurfaces if and only if
$F$ is a Morse family of functions.
\end{Pro}
\demo
By the assumption, there exists $\lambda (q,x,t)\in \mathcal{E}_{k+n+1}$ such that $\lambda (0)\not= 0$ and
\[
\mathcal{F}(q,x,t)=\lambda (q,x,t)(F(q,x)-t).
\]
Since $\partial \mathcal{F}/\partial q_i=\partial \lambda/\partial q_i(F-t)+\lambda \partial F/\partial q_i,$
we have
\[
\Delta ^*(\mathcal{F})=(\mathcal{F},d_1\mathcal{F})=\left(\lambda (F-t),\frac{\partial \lambda}{\partial q}(F-t)+\lambda \frac{\partial F}{\partial q}\right),
\]
where
\[
\frac{\partial \lambda}{\partial q}(F-t)+\lambda \frac{\partial F}{\partial q}=\left(
\frac{\partial \lambda}{\partial q_1}(F-t)+\lambda \frac{\partial F}{\partial q_1},\dots ,
\frac{\partial \lambda}{\partial q_k}(F-t)+\lambda \frac{\partial F}{\partial q_k}\right).
\]
By straightforward calculations, the Jacobian matrix of $\Delta ^*(\mathcal{F})(0)$ is
\[J_{\Delta ^*(\mathcal{F})}(0)=\begin{pmatrix}
0 &
\lambda (0) \frac{\partial F}{\partial x}(0) & -\lambda (0) \\
\lambda (0) \frac{\partial^2 F}{\partial q^2}(0) & \lambda (0) \frac{\partial^2 F}{\partial x\partial q}(0) &  0 
\end{pmatrix}
\]
We remark that the Jacobi matrix of $\Delta F$ is given by $J_{\Delta F}=(\partial^2 F/\partial q^2\  \partial^2 F/\partial x\partial q).$
Therefore, ${\rm rank}\, J_{\Delta ^*(\mathcal{F})}(0)=k+1$ if and only if ${\rm rank}\, J_{\Delta F}(0)= k.$
This completes the proof.
\enD
We now consider the case $\mathcal{F}(q,x,t)=\lambda (q,x,t)(F(q,x)-t).$ 
In this case, 
\[
\Sigma _*(\mathcal{F})=\{(q,x,F(q,x))\in (\R^k\times (\R^n\times\R),0)\ |\ (q,x)\in C(F)\},
\]
where
$
C(F)=\Delta F^{-1}(0).
$
Moreover, we have the Lagrangian submanifold germ $L(F)(C(F))\subset T^*\R^n$, where 
\[
L(F)(q,x)=\left(x,\frac{\partial F}{\partial x_1}(q,x),\dots ,\frac{\partial F}{\partial x_n}(q,x)\right).
\]
Since $\mathcal{F}$ is a graph-like Morse family of hypersurfaces, we have
a big Legendrian submanifold germ $\mathscr{L}_{\mathcal{F}}(\Sigma_{*}(\mathcal{F}))\subset J^1_{GA}(\R^n,\R),$ where
$\mathscr{L}_{\mathcal{F}}:(\Sigma_{*}(\mathcal{F}),0) \to J^1_{GA}(\R^n,\R)$ is defined by 
\[
\mathscr{L}_{\mathcal F}(q,x,t)=\left(x,t,-\frac{\displaystyle \frac{\partial \mathcal{F}}{\displaystyle\partial x_1}(q,x,t)}{\frac{\displaystyle \partial \mathcal{F}}{\displaystyle\partial t}(q,x,t)},\dots , -\frac{\displaystyle\frac{\partial \mathcal{F}}{\partial x_n}(q,x,t)}{\frac{\displaystyle\partial \mathcal{F}}{\displaystyle\partial t}(q,x,t)},\right)\in J^1_{GA}(\R^n,\R)\cong T^*\R^n\times \R.
\]
We also define a map ${\mathfrak L}_F:(C(F),0) \to J^1_{GA}(\R^n,\R)$ by
$$
{\mathfrak L}_F(q,x)=\left(x,F(q,x),\frac{\partial F}{\partial x_1}(q,x),\dots,\frac{\partial F}{\partial x_n}(q,x)\right).
$$
Since $\partial \mathcal{F}/\partial x_i=\partial \lambda/\partial x_i(F-t)+\lambda \partial F/\partial x_i$
and $\partial \mathcal{F}/\partial t=\partial \lambda/\partial t(F-t)-\lambda,$
we have $\partial \mathcal{F}/\partial x_i(q,x,t)= \lambda (q,x,t)\partial F/\partial x_i (q,x,t)$
and $\partial \mathcal{F}/\partial t(q,x,t)=-\lambda (q,x,t)$ for $(q,x,t)\in \Sigma _*(\mathcal{F}).$
It follows that
${\mathfrak L}_F(C(F))=\mathscr{L}_{\mathcal F}(\Sigma _*({\mathcal F})).$
By definition, we have $\Pi(\mathscr{L}_{\mathcal F}(\Sigma _*({\mathcal F})))=\Pi({\mathfrak L}_F(C(F)))=L(F)(C(F)).$
The graph-like wave front of the graph-like Legendrian unfolding ${\mathfrak L}_F(C(F))=\mathscr{L}_{\mathcal F}(\Sigma _*({\mathcal F}))$ is the graph of $F|_{C(F)}.$ 
This is the reason why we call it a graph-like Legendrian unfolding.
For a non-degenerate graph-like Morse family of hypersurfaces, we have the following proposition.
\begin{Pro}
With the same notations as Proposition 4.4, $\mathcal{F}$ is a non-degenerate graph-like Morse family of hypersurfaces
if and only if $F$ is a Morse family of hypersurfaces. In this case, $F$ is also a Morse family of functions such that
\[
\left(\frac{\partial F}{\partial x_1}(0),\dots ,\frac{\partial F}{\partial x_n}(0)\right)\not= \bo.
\]
\end{Pro}
\demo
By exactly the same calculations as those in the proof of Proposition 4.4, 
the Jacobi matrix of $\Delta ^*(\mathcal{F}|_{\R^k\times \R^n\times \{0\}})$ is
\[J_{\Delta ^*(\mathcal{F}|_{\R^k\times \R^n\times \{0\}})}(0)=\begin{pmatrix}
0 &
\lambda (0) \frac{\partial F}{\partial x}(0) ) \\
\lambda (0) \frac{\partial^2 F}{\partial q^2}(0) & \lambda (0) \frac{\partial^2 F}{\partial x\partial q}(0) 
\end{pmatrix}.
\]
On the other hand, the Jacobi matrix of $\Delta ^*(F)$ is 
\[J_{\Delta ^*(F)}(0)=\begin{pmatrix}
0 &
\frac{\partial F}{\partial x}(0) ) \\
\frac{\partial^2 F}{\partial q^2}(0) & \frac{\partial^2 F}{\partial x\partial q}(0) 
\end{pmatrix},
\]
so that the first assertion holds.
Moreover, 
${\rm rank}\, J_{\Delta ^*(\mathcal{F}|_{\R^k\times \R^n\times \{0\}})}(0)=k+1$ implies 
${\rm rank}\, J_{\Delta F}(0)= k$ and $\partial F/\partial x(0)\not= \bo.$
This completes the proof.
\enD
The momentary front for a fixed $t\in (\R,0)$ is 
$
W_t(L)=\pi_1 (\pi_2^{-1}(t)\cap W(L))$.
We define $L_t=L\cap (\pi_2\circ\overline{\pi})^{-1}(t)= L\cap (T^*\R^n\times \{t\})$ under the canonical identification $J^1_{GA}(\R^n,\R)\cong T^*\R^n\times \R.$
Then $\Pi(L)\subset T^*\R^n$ and $\widetilde{\pi}\circ\Pi (L_t)\subset PT^*\R^n,$ where
$\widetilde{\pi}:T^*\R^n\lon PT^*(\R^n)$ is the canonical projection.
We also have the canonical projections $\varpi:T^*\R^n\lon \R^n$ and $\overline{\varpi}:PT^*\R^n\lon \R^n$
such that $\pi _1\circ \overline{\pi}=\varpi\circ \Pi$ and $\overline{\varpi}\circ\widetilde{\pi}=\varpi.$
Then we have the following proposition.
\begin{Pro} Let $L\subset J^1_{GA}(\R^n,\R)$ be a non-degenerate graph-like Legendrian unfolding.
Then $\Pi (L)$ is a Lagrangian submanifold and $\widetilde{\pi}\circ \Pi(L_t)$ is a Legendrian submanifold
in $PT^*(\R^n).$
\end{Pro}
\demo
 By Proposition 4.1, $\Pi (L)$ is a Lagrangian submanifold in $T^*\R^n.$
 Since $L$ is a non-degenerate Legendrian unfolding in $J^1_{GA}(\R^n,\R),$ we have a non-degenerate graph-like generating family $\mathcal{F}$ of $L$ at least locally. This means that $L=\mathscr{L}_{\mathcal{F}}(\Sigma _*(\mathcal{F}))$ as set germs.
 Since $\mathcal{F}$ is a graph-like Morse family of hypersurface, it is written as 
 $\mathcal{F}(q,x,t)=\lambda (q,x,t)(F(q,x)-t).$
Therefore, we have $\mathscr{L}_{\mathcal F}(\Sigma _*({\mathcal F}))={\mathfrak L}_F(C(F)).$
By definition, $\Pi\circ {\mathfrak L}_F(C(F))=L(F)(C(F)),$ so that $F$ is a generating family of
$\Pi (L)$, locally. By Proposition 4.5, $F$ is also a Morse family of hypersurface, so
that $\mathscr{L}_F(\Sigma _*(F))$ is a Legendrian submanifold germ in $PT^*(\R^n).$
Without loss of generality, we can assume that $t=0.$
Since $\Sigma _*(F)=C(F)\cap F^{-1}(0),$ $\mathscr{L}_F(\Sigma _*(F))=\widetilde{\pi}\circ\Pi (\mathfrak{L}_F(C(F))\cap (\pi_2\circ\overline{\pi})^{-1}(0))=\widetilde{\pi}\circ\Pi(L_0).$
This completes the proof.
\enD
In general, the momentary front $W_t(L)$ of a big Legendrian submanifold $L\subset PT^*(\R^n\times\R)$ is not necessarily a wave front of a Legendrian submanifold in the ordinary sense.
However, for a non-degenerate Legendrian unfolding in $J^1_{GA}(\R^n,\R),$ we have the following corollary.
\begin{Co} Let $L\subset J^1_{GA}(\R^n,\R)$ be a non-degenerate graph-like Legendrian unfolding.
Then the momentary front $W_t(L)$ is the wave front set of the Legendrian submanifold $\widetilde{\pi}\circ\Pi(L_t)\subset PT^*(\R^n).$
Moreover, the caustic $C_L$ is the caustic of the Lagrangian submanifold $\Pi (L)\subset T^*\R^n.$
In other words, $W_t(L)=\overline{\varpi}(\widetilde{\pi}\circ\Pi(L_t))$ and $C_L$ is the singular value set of $\varpi |\Pi(L).$
\end{Co}
\demo
By definition, we have
\[
\overline{\pi}(L_t)=\overline{\pi}(L\cap (\pi_2\circ\overline{\pi})^{-1}(t))=W(L)\cap \pi_2^{-1}(t),
\]
so that 
\[
W_t(L)=\pi _1(W(L)\cap \pi_2^{-1}(t))=\pi_1\circ\overline{\pi}(L_t)=\varpi\circ\Pi(L_t)=\overline{\varpi}(\widetilde{\pi}\circ\Pi(L_t)).
\]
\par
We remark that $\pi_1\circ\overline{\pi}=\varpi\circ\Pi$. 
By Proposition 4.1, $z\in L$ is a singular point of $\overline{\pi}|_L:L\lon \R^n\times\R$
if and only if it is a singular point of $\varpi|_{\Pi(L)}:\Pi(L)\lon \R^n.$
Therefore, the caustic $C_L$ is the singular value set of $\varpi |_{\Pi(L)}.$
\enD
For a graph-like Morse family of hypersurfaces $\mathcal{F}(q,x,t)=\lambda (q,x,t)(F(q,x)-t),$
$\mathcal{F}(q,x,t)$ and $\overline{F}(q,x,t)=F(q,x)-t$ are $s$-$S.P^+$-$\mathcal{K}$-equivalent,
so that we consider $\overline{F}(q,x,t)=F(q,x)-t$ as a graph-like Morse family.
Moreover, if $\mathcal{F}$ is non-degenerate, then
$F(q,x)$ is a Morse family of functions.
We now suppose that $F(q,x)$ is a Morse family of functions.
Consider the graph-like Morse family of hypersurfaces $\overline{F}(q,x,t)=F(q,x)-t$
which is not necessarily non-degenerate.
Then we have $\mathscr{L}_{\overline{F}}(\Sigma _*(\overline{F}))=\mathfrak{L}_F(C(F)).$
We also denote that $\overline{f}(q,t)=f(q)-t$ for any $f\in \mathfrak{M}_k.$ 
We can represent the extended tangent space of $\overline{f}:(\R^k\times \R,0)\lon (\R,0)$
 relative to $S.P^+$-${\cal K}$ by
$$
T_e(S.P^+\mbox{\rm -}{\cal K})(\overline{f})=\left\langle \frac{\partial f}{\partial q_1}(q),\dots ,
\frac{\partial f}{\partial q_k}(q) ,f(q)-t\right\rangle _{{\cal E}_{(q,t)}}+
\langle1\rangle _{\R}.
$$
For a deformation $\overline{F}:(\R^k\times\R^n \times \R,0)\lon (\R,0)$ of $\overline{f},$ $\overline{F}$ is 
infinitesimally $S.P^+$-${\cal K}$-versal deformation of $\overline{f}$ if and only if 
$$
{\cal E}_{(q,t)}=T_e(S.P^+\mbox{\rm -}{\cal K})(\overline{f})+
\left\langle \frac{\partial F}{\partial x_1}|_{\R^k\times\{0\}},\dots
,\frac{\partial F}{\partial x_n}|_{\R^k\times\{0\}}\right\rangle _{\R}.
$$

We compare the equivalence relations 
between Lagrangian submanifold germs and induced graph-like Legendrian unfoldings, that is, between Morse families of functions and graph-like Morse families of hypersurfaces. 
As a consequence, we give a relationship between caustics and graph-like wave fronts. 

\begin{Pro}[\cite{Izumiya-Takahashi}]
If Lagrangian submanifold germs $L(F)(C(F))$, $L(G)(C(G))$ are Lagrangian equivalent,
then the graph-like Legendrian unfoldings ${\mathfrak L}_F(C(F))$, ${\mathfrak L}_G(C(G))$ are $S.P^+$-Legendrian equivalent.
\end{Pro}

\demo By Proposition 2.1, two Lagrangian submanifold germs $L(F)(C(F))$, $L(G)(C(G))$
are Lagrangian equivalent if and only if $F$ and $G$ are stably
$P$-${\mathcal R}^+$-equivalent.
By definition, if $F$ and $G$ are stably $P$-${\mathcal R}^+$-equivalent,
then $\overline{F}$ and $\overline{G}$ are stably $s$-$S.P^+$-${\cal K}$-equivalent.
 By the assertion (1) of Theorem 3.5, ${\mathfrak L}_F(C(F))$ and ${\mathfrak L}_G(C(G))$ are  $S.P^+$-Legendrian equivalent.
\enD

\begin{Rem}{\rm 
The above proposition asserts that Lagrangian equivalence is a stronger equivalence relation than $S.P^+$-Legendrian equivalence.
The $S.P^+$-Legendrian equivalence relation among graph-like Legendrian unfoldings preserves both the diffeomorphism types of  caustics and Maxwell stratified sets.
On the other hand, if we observe the real caustics of rays, we cannot observe the structure of wave front propagations and the Maxwell stratified sets.
In this sense, there are hidden structures behind the picture of real caustics. 
By the above proposition, Lagrangian equivalence preserves not only the diffeomorphism type of caustics, but also the hidden geometric structure of wave front propagations. }
\end{Rem}

It seems that the converse assertion does not hold.
However, we have the following proposition.
\begin{Pro}[\cite{Izumiya-Takahashi2}]
Suppose that $L(F)(C(F))$ and $L(G)(C(G))$ are Lagrange stable. 
If the graph-like Legendrian unfoldings ${\mathfrak L}_F(C(F))$ and ${\mathfrak L}_G(C(G))$ are $S.P^+$-Legendrian equivalent, then the Lagrangian submanifold germs $L(F)(C(F))$ and  
$L(G)(C(G))
$ are Lagrangian equivalent.
\end{Pro}

In order to prove the proposition, we need the following lemma:
\begin{Lem}
If $\overline{f}$ and $\overline{g}:(\R^k \times \R,0) \lon (\R,0)$ are $S.P$-$\mathcal{K}$-equivalent, then $f$ and $g:(\R^k,0) \lon (\R,0)$ are $\mathcal{R}$-equivalent, where $\overline{f}(q,t)=f(q)-t$ and $\overline{g}(q,t)=g(q)-t$. 
\end{Lem}

\demo
By the definition of $S.P$-$\mathcal{K}$-equivalence, there exist a diffeomorphism germ of $\Phi:(\R^k \times \R,0) \lon (\R^k \times \R,0)$ of the form $\Phi(q,t)=(\phi(q,t),t)$ and a non-zero function germ $\lambda:(\R^k \times \R,0) \lon \R$ such that $\overline{f}= \lambda \cdot \overline{g} \circ \Phi$. 
Then the diffeomorphism $\Phi$ preserves the zero-level set of $\overline{f}$ and $\overline{g}$, that is, $\Phi (\overline{f}^{-1}(0))=\overline{g}^{-1}(0)$. 
Since the zero-level set of $\overline{f}$ is the graph of $f$ and the form of $\Phi$, we have $f=g \circ \psi$, where $\psi(q)=\phi(q,f(q))$. 
It is easy to show that $\psi:(\R^k,0) \lon (\R^k,0)$ is a diffeomorphism germ. 
Hence $f$ and $g$ are $\mathcal{R}$-equivalent.
\enD

{\it Proof of Proposition 4.10}. 
By the assertion (1) of Theorem 3.5, $\overline{F}$ and $\overline{G}$ are stably $s$-$S.P^+$-${\cal K}$-equivalent. 
It follows that $\overline{f}$ and $\overline{g}$ are stably $S.P$-$\mathcal{K}$-equivalent. 
By Lemma 4.11, $f$ and $g$ are stably $\mathcal{R}$-equivalent. 
By the uniqueness of the infinitesimally $\mathcal{R}^+$-versal unfolding (cf., \cite{Bro}), $F$ and $G$ are stably $P$-$\mathcal{R}^+$-equivalent.
\enD

By definition, the set of Legendrian singular points of a graph-like Legendrian unfolding $\mathfrak{L}_F(C(F))$ coincides with the set of singular points of $\pi \circ L(F)$. 
Therefore the singularities of graph-like wave fronts of $\mathfrak{L}_F(C(F))$ lie on the caustics of $L(F)$. 
Moreover, if Lagrangian submanifold germ $L(F)(C(F))$ is Lagrangian stable, 
then the regular set of $\overline{\pi} \circ \mathfrak{L}_F(C(F))$ is dense. Hence we can apply Proposition 3.1 to our situation and obtain the following theorem as a corollary of Propositions 4.8 and 4.10. 

\begin{Th}[\cite{Izumiya-Takahashi2}]
Suppose that $L(F)(C(F))$ and $L(G)(C(G))$ are Lagrangian stable. 
Then Lagrangian submanifold germs $L(F)(C(F))$ and $L(G)(C(G))$ are Lagrangian equivalent if and only if graph-like wave fronts $W(\mathfrak{L}_F)$ and $W(\mathfrak{L}_G)$ are $S.P^+$-diffeomorphic.
\end{Th}
\par
Moreover, we have the following theorem.
\begin{Th}[\cite{Izumiya-Takahashi3}]
Suppose that $\mathcal{F}(q,x,t)=\lambda(q,x,t)\langle F(q,x)-t\rangle$ is a graph-like Morse family of hypersurfaces. Then
$\mathscr{L}_{\mathcal{F}}(\Sigma _*(\mathcal{F}))$ is $S.P^+$-Legendrian stable if and only if $L(F)(C(F))$ is Lagrangian stable.
\end{Th}
\demo
By Proposition 4.8, if $L(F)(C(F))$ is Lagrangian stable, then $\mathscr{L}_{\mathcal{F}}(\Sigma _*(\mathcal{F}))$ is $S.P^+$-Legendrian stable.
For the converse, suppose that ${\mathfrak L}_F(C(F))$ is a $S.P^+$-Legendre stable.
By the assertion (2) of Theorem 3.5, we have
$$
{\rm dim}_\R \frac{\mathcal{E}_{k+1}}{\langle \frac{\partial f}{\partial q_1}(q),\dots ,
\frac{\partial f}{\partial q_k}(q) ,f(q)-t \rangle_{\mathcal{E}_{k+1}}+\langle1\rangle _{\R}}< \infty.
$$
It follows that ${\rm dim}_{\R} \mathcal{E}_k/\langle \frac{\partial f}{\partial q_1}(q),\dots ,
\frac{\partial f}{\partial q_k}(q) ,f(q)\rangle_{\mathcal{E}_{k}}<\infty$, namely, 
$f$ is a $\mathcal{K}$-finitely determined (see the definition \cite{Dimca,Martinet}).
It is a well-known that $f$ is $\mathcal{K}$-finitely determined if and only if $f$ is $\mathcal{R}^+$-finitely determined, see \cite{Dimca}.
Under the condition that $f$ is $\mathcal{R}^+$-finitely determined, 
$F$ is an infinitesimally $\mathcal{R}^+$-versal deformation of $f$ if and only if $F$ is an $\mathcal{R}^+$-transversal deformation of $f$, namely, there exists a number $\ell \in \mathbb{N}$ such that
\begin{eqnarray}\label{transversal}
\mathcal{E}_k=J_f + {\left\langle 
\frac{\partial{F}}{\partial{x_1}} | _{\R^k \times \{0\}},\dots,
\frac{\partial{F}}{\partial{x_n}} | _{\R^k \times \{0\}}
\right\rangle}_\R +{\langle 1 \rangle}_\R+\mathcal{M}^{\ell+1}_k.
\end{eqnarray}
Hence, it is enough to show the equality (\ref{transversal}).
Let $g(q) \in \mathcal{E}_k$. 
Since $g(q) \in \mathcal{E}_{k+1}$, there exist $\lambda_i(q,t), \mu(q,t) \in \mathcal{E}_{k+1}$ $(i=1,\dots,k)$ and $c, c_j \in \R$ $(j=1,\dots,n)$ such that 
\begin{eqnarray}\label{eq}
g(q)=\sum^{k}_{i=1} \lambda_i(q,t)\frac{\partial f}{\partial q_i}(q)+\mu(q,t)(f(q)-t)+c+\sum^n_{j=1} c_j \frac{\partial F}{\partial x_j}(q,0).
\end{eqnarray}
Differentiating the equality (\ref{eq}) with respect to $t$, we have
\begin{eqnarray}\label{diff.t}
0=\sum^k_{i=1} \frac{\partial \lambda_i}{\partial t} (q,t) \frac{\partial f}{\partial q_i} (q)+\frac{\partial \mu}{\partial t}(q,t)(f(q)-t)-\mu (q,t).
\end{eqnarray}
We put $t=0$ in (\ref{diff.t}), $0=\sum^k_{i=1} ({\partial \lambda_i}/{\partial t}) (q,0) ({\partial f}/{\partial q_i}) (q)+({\partial \mu}/{\partial t})(q,0)f(q)-\mu (q,0)$.
Also we put $t=0$ in (\ref{eq}), then 
\begin{eqnarray}\label{eq2}
g(q) &=&  \sum^{k}_{i=1} \lambda_i(q,0)\frac{\partial f}{\partial q_i}(q)+\mu(q,0)f(q)+c+\sum^n_{j=1} c_j \frac{\partial F}{\partial x_j}(q,0)\nonumber \\
&=& \sum^k_{i=1} \alpha_i (q) \frac{\partial f}{\partial q_i}(q)+\frac{\partial \mu}{\partial t}(q,0)f^2(q)+c+\sum^n_{j=1} c_j \frac{\partial F}{\partial x_j}(q,0),
\end{eqnarray}
for some $\alpha_i \in \mathcal{E}_k, i =1\dots,k$.
Again differentiating (\ref{diff.t}) with respect to $t$ and put $t=0$, then
$$
0=\sum^k_{i=1} \frac{\partial^2 \lambda_i}{\partial t^2} (q,0) \frac{\partial f}{\partial q_i} (q)+\frac{\partial^2 \mu}{\partial t^2}(q,0)f(q)-2\frac{\partial \mu}{\partial t} (q,0).
$$
Hence (\ref{eq2}) is equal to 
$$
\sum^k_{i=1} \beta_i (q) \frac{\partial f}{\partial q_i}(q)+\frac{1}{2} \frac{\partial^2 \mu}{\partial t^2}(q,0)f^3(q)+c+\sum^n_{j=1} c_j \frac{\partial F}{\partial x_j}(q,0),
$$
for some $\beta_i \in \mathcal{E}_k, i=1,\dots,k$. 
Inductively, we take $\ell$-times differentiate (\ref{diff.t}) with respect to $t$ and put $t=0$, then we have 
$$
g(q)=\sum^k_{i=1} \gamma_i (q) \frac{\partial f}{\partial q_i}(q)+\frac{1}{\ell!}\frac{\partial^{\ell} \mu}{\partial t^{\ell}}(q,0)f^{\ell+1}(q)+c+\sum^n_{j=1} c_j \frac{\partial F}{\partial x_j}(q,0),
$$
for some $\gamma_i \in \mathcal{E}_k, i=1,\dots,k$. 
It follows that $g(q)$ is contained in the right hand of (\ref{transversal}).
This completes the proof.
\enD
\par
One of the consequences of the above arguments is the following theorem on the relation among
graph-like Legendrian unfoldings and Lagrangian singularities.

\begin{Th} Let $\mathcal{F}:(\R^k\times \R^n\times \R, 0)\lon (\R,0)$ and
$\mathcal{G}:(\R^{k'}\times \R^n\times \R,0)\lon (\R,0)$ be graph-like Morse families of hypersurfaces of the forms $\mathcal{F}(q,x,t)=\lambda(q,x,t)(F(q,x)-t)$ and
$\mathcal{G}(q',x,t)=\mu(q',x,t)(G(q',x)-t)$ such that
$\mathscr{L}_{\mathcal{F}}(\Sigma _*(\mathcal{F}))$
and $\mathscr{L}_{\mathcal{G}}(\Sigma _*(\mathcal{G}))$ are $S.P^+$-Legendrian stable. Then the following conditions are
equivalent{\rm :}
\par\noindent
{\rm (1)} $\mathscr{L}_{\mathcal{F}}(\Sigma _*(\mathcal{F}))$
and $\mathscr{L}_{\mathcal{G}}(\Sigma _*(\mathcal{G}))$ are $S.P^+$-Legendrian equivalent,
\par\noindent
{\rm (2)} $\mathcal{F}$ and $\mathcal{G}$ are stably $s$-$S.P^+$-$\mathcal{K}$-equivalent,
\par\noindent
{\rm (3)} $\overline{f}(q,t)=F(q,0)-t$ and $\overline{g}(q',t)=G(q',0)-t$ are stably $S.P$-$\mathcal{K}$-equivalent,
\par\noindent
{\rm (4)} $f(q)=F(q,0)$ and $g(q')=G(q',0)$ are stably $\mathcal{R}$-equivalent,
\par\noindent
{\rm (5)} $F(q,x)$ and $G(q',x)$ are stably $P$-$\mathcal{R}^+$-equivalent,
\par\noindent
{\rm (6)} $L(F)(C(F))$ and $L(G)(C(G))$ are Lagrangian equivalent,
\par\noindent
{\rm (7)} $W(\mathscr{L}_{\mathcal{F}}(\Sigma _*(\mathcal{F})))$ and $W(\mathscr{L}_{\mathcal{G}}(\Sigma _*(\mathcal{G})))$
are $S.P^+$-diffeomorphic.
\end{Th}
\demo
By the assertion (1) of Theorem 3.5,  the conditions (1) and (2) are equivalent.
By definition, the condition (2) implies the condition (3), the condition (4) implies (3)
and the condition (5) implies (4), respectively.
By Lemma 4.11, the condition (3) implies the condition (4).
By Theorem 2.1, the conditions (5) and (6) are equivalent.
It also follows from the definition that the condition (1) implies (7).
We remark that all these assertions hold without the assumptions of the $S.P^+$-Legendrian stability.
Generically, the condition (7) implies the condition (1) by Proposition 3.1.
Of course, the assertion of Theorem 4.12 holds under the assumption of $S.P^+$-Legendrian stability.
By the assumption of $S.P^+$-Legendrian stability, the graph-like Morse families of
hypersurface $\mathcal{F}$ and $\mathcal{G}$ are infinitesimally $S.P^+$-$\mathcal{K}$-versal deformations
of  $\overline{f}$ and  $\overline{g}$, respectively (cf., Theorem 3.5, (2)).
By the uniqueness result for infinitesimally $S.P^+$-$\mathcal{K}$-versal deformations, the condition (3) implies the condition (2).
Moreover, by Theorem 4.13, $L(F)(C(F))$ and $L(G)(C(G))$ are Lagrangian stable.
This means that $F$ and $G$ are infinitesimally $\mathcal{R}^+$-versal deformations of $f$ and $g$, respectively. Therefore by the uniqueness results for  infinitesimally $\mathcal{R}^+$-versal deformations,
the condition (4) implies the condition (5).
This completes the proof.
\enD
\begin{Rem}{\rm 
(1) By Theorem 4.13, the assumption of the above theorem is equivalent to the condition that
$L(F)(C(F))$ and $L(G)(C(G))$ are Lagrangian stable.
\par\noindent
(2) If $k=k'$ and $q=q'$ in the above theorem, we can remove the word \lq\lq stably\rq\rq
\ in the conditions (2),(3),(4) and (5).
\par\noindent
(3) The $S.P^+$-Legendrian stability of $\mathscr{L}_{\mathcal{F}}(\Sigma _*(\mathcal{F}))$ is a generic condition for $n\leq 5.$ 
\par\noindent
(4) By the remark in the proof of the above theorem, the conditions (1) and (7) are equivalent
generically for an arbitrary dimension $n$ without the assumption on the $S.P^+$-Legendrian stability.
Therefore, the conditions (1),(2) and (7) are all equivalent to each other as before.
Lagrangian equivalence (i.e., the conditions (5) and (6)) is a stronger condition than others
as before.
}
\end{Rem}

\section{Applications}
In this section we explain some applications of the theory of wave front propagations.
\subsection{Completely integrable first order ordinary differential equations}
In this subsection we consider implicit first order ordinary differential equations.
There are classically written as $F(x,y,dy/dx)=0.$
However, if we set $p=dy/dx,$ then we have a surface in the 1-jet space $J^1(\R,\R)$
defined by $F(x,y,p)=0,$ where we have the canonical contact form $\theta =dy-pdx.$
Generically, we may assume that the surface is regular, then it has a local parametrization, so that
it is the image of an immersion
at least locally.
An {\it ordinary differential equation germ} (briefly, an {\it ODE}) is defined to be
an immersion germ $f:(\R^2,0)\lon J^1(\R,\R).$
We say that an ODE $f:(\R^2,0)\lon J^1(\R,\R)$ is {\it completely integrable}
if there exists a submersion germ $\mu :(\R^2,0)\lon (\R,0)$ such that
$f^*\theta \in \langle d\mu\rangle _{\mathcal{E}_2}.$
It follows that there exists a unique $h\in \mathcal{E}_2$ such that
$f^*\theta =hd\mu.$
In this case we call $\mu$ a {\it complete integral} of $f.$
In \cite{HIIY} a generic classification has been considered of completely integrable first order ODEs 
by point transformations.
Let $f,g:(\R^2,0)\lon J^1(\R,\R)\subset PT^*(\R\times \R)$ be ODEs.
We say that $f,g$ are {\it equivalent} as ODEs if there exist diffeomorphism germs $\psi :(\R^2,0)\lon (\R^2,0)$
and $\Phi :(\R\times \R,0)\lon (\R\times \R,0)$ such that $\widehat{\Phi}\circ f=g\circ \psi.$
Here $\widehat{\Phi}$ is the unique contact lift of $\Phi.$
The diffeomorphism germ $\Phi :(\R\times\R,0)\lon (\R\times \R,0)$ is traditionally called 
a {\it point transformation}.
We represent $f$ by the canonical coordinates of $J^1(\R,\R)$ by
$f(u_1,u_2)=(x(u_1,u_2),y(u_1,u_2),p(u_1,u_2)).$ 
If we have a complete integral $\mu:(\R^2,0)\lon \R$ of $f$, we 
define a immersion germ $\ell _{(\mu,f)}:(\R^2,0)\lon J^1(\R\times\R,\R)$
by $$\ell_{(\mu,f)}(u_1,u_2)=(\mu (u_1,u_2),x(u_1,u_2),y(u_1,u_2),h(u_1,u_2),p(u_1,u_2)).$$
Then we have $\ell_{(\mu,f)}^*\Theta =0,$ for $\Theta =dy-pdx-qdt,$ where $(t,x,y,q,p)$ is
the canonical coordinate of $J^1(\R\times\R,\R).$
Therefore, the image of $\ell_{(\mu,f)}$ is a big Legendrian submanifold germ of $J^1(\R\times\R,\R).$
Here, we consider the parameter $t$ as the time-parameter.
Since the contact structure is defined by the contact form $\Theta =dy-pdx-qdt,$ 
$J^1(\R\times\R,\R)$ is of course an affine coordinate neighbourhood of $PT^*(\R\times\R\times\R)$ but it is not equal to $J^1_{AG}(\R\times\R,\R)\subset PT^*(\R\times\R\times \R).$
The above notation induces a divergent diagram of map germs as follows:
\[
\R\overset {\pi _1\circ\overline{\pi}\circ \ell_{(\mu,f)}}\longleftarrow  (\R^2,0)\overset {\pi _2\circ\overline{\pi}\circ \ell_{(\mu,f)}}\lon (\R\times\R,0),
\]
where $\overline{\pi} :J^1(\R\times \R,\R)\lon \R\times\R\times \R$ is $\overline{\pi}(t,x,y,q,p)=(t,x,y),$
$\pi _1:(\R\times\R\times\R,0)\lon \R$ is $\pi _1(t,x,y)=t$ and $\pi _2:(\R\times\R\times\R,0)\lon \R\times\R$ is $\pi _2(t,x,y)=(x,y).$
Actually, we have $\pi _1\circ\overline{\pi}\circ \ell_{(\mu,f)}=\mu$ and $\pi _2\circ\overline{\pi}\circ \ell_{(\mu,f)}=\widehat{\pi}\circ f$,
where $\widehat{\pi}:J^1(\R,\R)\lon \R\times\R$ is the canonical projection $\widehat{\pi}(x,y,p)=(x,y).$
The space of completely integrable ODEs is identified with the space of big Legendrian submanifold
such that the restrictions of the $\pi _1\circ \overline{\pi}$-projection are non-singular.
For a divergent diagram
\[
\R\overset {\mu}\longleftarrow  (\R^2,0)\overset {g}\lon (\R\times\R,0),
\]
we say that $(\mu, g)$ is an {\it integral diagram} if there exist an immersion germ $f:(\R^2,0)\lon J^1(\R,\R)$
and a submersion germ $\mu :(\R^2,0)\lon \R$ such that
$
g=\widehat{\pi}\circ f.
$
Therefore we can apply the theory of big wave fronts.
In \cite{HIIY} the following proposition has been shown.
\begin{Pro}[\cite{HIIY}]
Let $f_i$ $(i=1,2)$ be completely integrable first order ODEs with the integrals $\mu _i$ and
the corresponding integral diagrams are $(\mu_i,g_i)$. Suppose that sets of Legendrian singular points
of $\ell_{(\mu_i,f_i)}$ $(i=1,2)$ are nowhere dense. Then the following conditions are equivalent:
\par\noindent
{\rm (1)} $f_1,f_2$ are equivalent as ODEs.
\par\noindent
{\rm (2)} 
There exists a diffeomorphism germ
$\Phi :(\R\times\R\times \R,0)\lon (\R\times\R\times \R,0)$ of the form $\Phi (t,x,y)=(\phi_1(t),\phi _2(x,y),\phi _3(x,y))$ such that $\widehat{\Phi}({\rm Image}\,\ell _{(\mu_1,f_1)})={\rm Image}\,\ell _{(\mu_2,f_2)}.$
\par\noindent
{\rm (3)} There exist diffeomorphism germs $\phi:(\R,0)\lon (\R,0)$, $\Phi: (\R^2,0)\lon (\R^2,0)$ and $\Psi :(\R\times\R,0)\lon (\R\times\R,0)$
such that $\phi\circ \mu_1 =\mu_2 \circ \Phi$
and $\Psi \circ g_1=g_2 \circ \Phi.$
\end{Pro}
\par
We say that two integral diagrams $(\mu_1,g_1)$ and $(\mu _2,g_2)$ are {\it equivalent} as integral diagrams
if the condition (3) of the above theorem holds.
By Remark 3.7, the classification by the above equivalence is almost impossible.
We also say that integral diagrams $(\mu_1,g_1)$ and $(\mu _2,g_2)$ are {\it strictly equivalent}
if the condition (3) of the above theorem holds for $\phi =1_{\R}$.
The strict equivalence corresponds to the $S.P$-Legendrian equivalence among the big Legendrian submanifold germs $\ell _{(\mu,f)}.$
Instead of the above equivalence relation, $S.P$-Legendrian equivalence was used for classification
in \cite{HIIY}. The technique used there was very hard.
In \cite{Izu95}, $S.P^+$-Legendrian equivalence was used.
If we have a classification of $\ell _{(\mu,f)}$ under the $S.P^+$-Legendrian equivalence,
we can automatically obtain the classification of integral diagrams by the strict equivalence.
\begin{Th}[\cite{HIIY,Izu95}]
For a \lq\lq generic\rq\rq  first order ODE $f:(\R^2,0)\lon J^1(\R,\R)$ with a complete integral $\mu:(\R^2,0)\lon \R,$
the corresponding integral diagram $(\mu,g)$ is strictly equivalent to one of the germs in the following list:
\par
{\rm (1)} $\mu=u_2$, $g=(u_1,u_2),$
\par
{\rm (2)} $\displaystyle{\mu=\frac{2}{3}u_1^3+u_2}$, $g=(u^2_1,u_2),$
\par
{\rm (3)} $\displaystyle{\mu=u_2-\frac{1}{2}u_1}$, $g=(u_1,u^2_2),$
\par
{\rm (4)} $\displaystyle{\mu=\frac{3}{4}u_1^4+\frac{1}{2}u_1^2u_2+u_2+\alpha\circ g}$, $g=(u^3_1+u_2u_1,u_2),$
\par
{\rm (5)} $\displaystyle{\mu=u_2+\alpha\circ g}$, $g=(u_1,u_2^3+u_1u_2),$
\par
{\rm (6)} $\displaystyle{\mu=-3u_2^2+4u_1u_2+u_1+\alpha\circ g}$, $g=(u_1,u^3_2+u_1u_2^2).$
\par\noindent
Here, $\alpha(v_1,v_2)$ are $C^\infty$-function germs, which are called functional modulus.
\end{Th}
\begin{Rem}{\rm The results has been generalized into the case for completely integrable holonomic systems
of first order partial differential equations \cite{Izu95,Izu-Kro}.
}
\end{Rem}
\par
In the list of the above theorem, the normal forms (3), (5) are said to be of {\it Clairaut type}.
The complete solutions for those equations are non-singular and the singular solutions are the envelopes of
the graph of complete solutions.
We say that a complete integrable first order ODE  $f:(\R^2,0)\lon J^1(\R,\R)$ with an integral $\mu:(\R^2,0)\lon \R$
is {\it Clairaut type} if $\widehat{\pi}\circ f|_{\mu^{-1}(t)}$ is non-singular for any $t\in \R$.
Then $\overline{\pi} \circ \ell _{(\mu,f)}$ is also non-singular.
In this case the discriminant of the family $\{W_t(\ell _{(\mu,f)}(\R^2))\}_{t\in (\R,0)}$ is equal to the envelope of the family of momentary
fronts $\Delta_{\ell _{(\mu,f)}(\R^2)}.$ Here, the momentary front is a special solution of the complete solution $\{\widehat{\pi}\circ f(\mu^{-1}(t))\}_{t\in \R}$.
This means that $\ell _{(\mu,f)}(\R^2)\cap J^1_{GA}(\R\times\R,\R)=\emptyset.$
\par
On the other hand, the normal forms (2), (4) are said to be of {\it regular type}.
In those cases $f^*\theta\not= 0$ and we have $\ell _{(\mu,f)}(\R^2)\subset  J^1_{GA}(\R\times\R,\R).$ 
Therefore, $\ell _{(\mu,f)}(\R^2)$ is a graph-like Legendrian unfolding, so that
the discriminant of the family $\{W_t(\ell _{(\mu,f)}(\R^2))\}_{t\in (\R,0)}$ is $C_{\ell _{(\mu,f)}(\R^2)}\cup M_{\ell _{(\mu,f)}(\R^2)}.$
Finally the normal form (6) is as before a {\it mixed hold type}.
In this case, $\ell _{(\mu,f)}(\R^2)\subset  J^1(\R\times\R,\R)$ but $\ell _{(\mu,f)}(\R^2)\not\subset  J^1_{GA}(\R\times\R,\R)$.
Actually, $\ell _{(\mu,f)}(0)\in \overline{J^1_{GA}(\R\times\R,\R)}$, where $\overline{X}$ is the closure of $X$.
The pictures of the families of momentary fronts of (4), (5), (6) are drawn in Figures 5, 6 and 7.
We can observe that the discriminants of the families $\{W_t(\ell _{(\mu,f)}(\R^2))\}_{t\in (\R,0)}$ are $C_{\ell _{(\mu,f)}(\R^2)}\cup M_{\ell _{(\mu,f)}(\R^2)}$ for (4), $\Delta_{\ell _{(\mu,f)}(\R^2)}$
for (5) and $C_{\ell _{(\mu,f)}(\R^2)}\cup \Delta_{\ell _{(\mu,f)}(\R^2)}$ for (6), respectively.
Moreover, the $C_{\ell _{(\mu,f)}(\R^2)}$ of the germ (4) and  $\Delta _{\ell _{(\mu,f)}(\R^2)}$ of the germ (5) are semi-cubical parabolas.
Therefore, these are diffeomorphic but their discriminants are not $S.P^+$-diffeomorphic.
\begin{figure}[htbp]
  \begin{center}
    \small
    \begin{tabular}{ccc}
       \includegraphics[width=40mm]{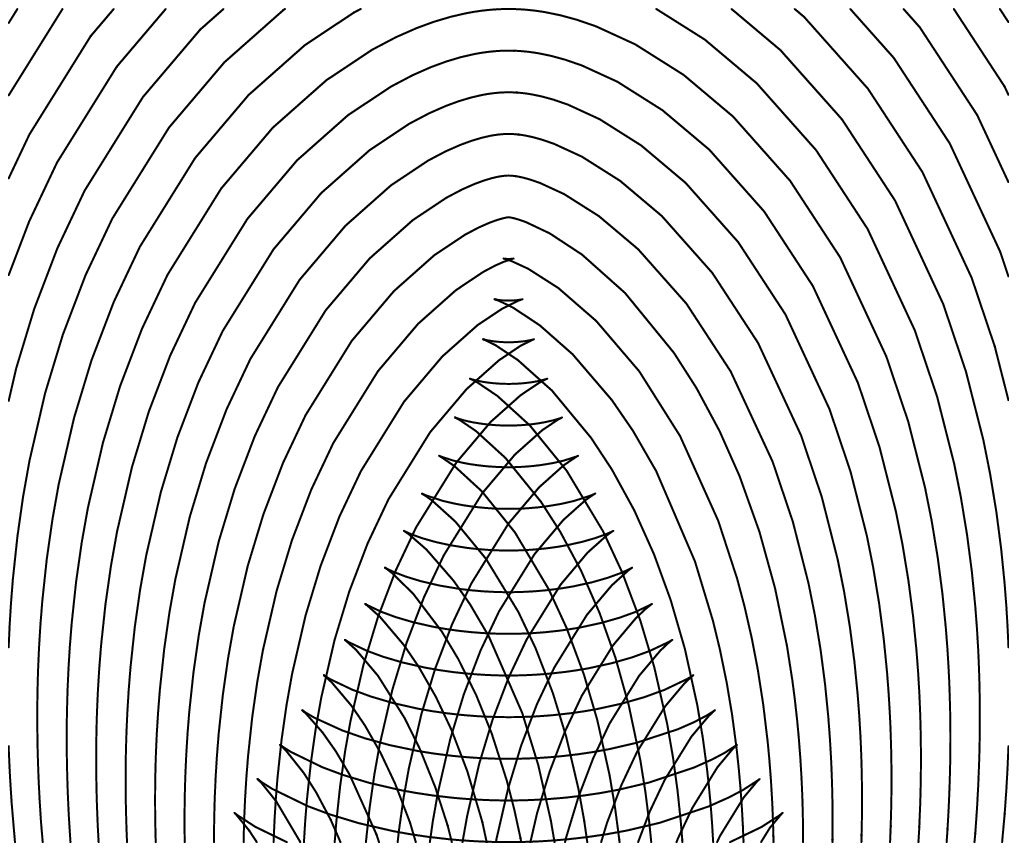} &
      \includegraphics[width=40mm]{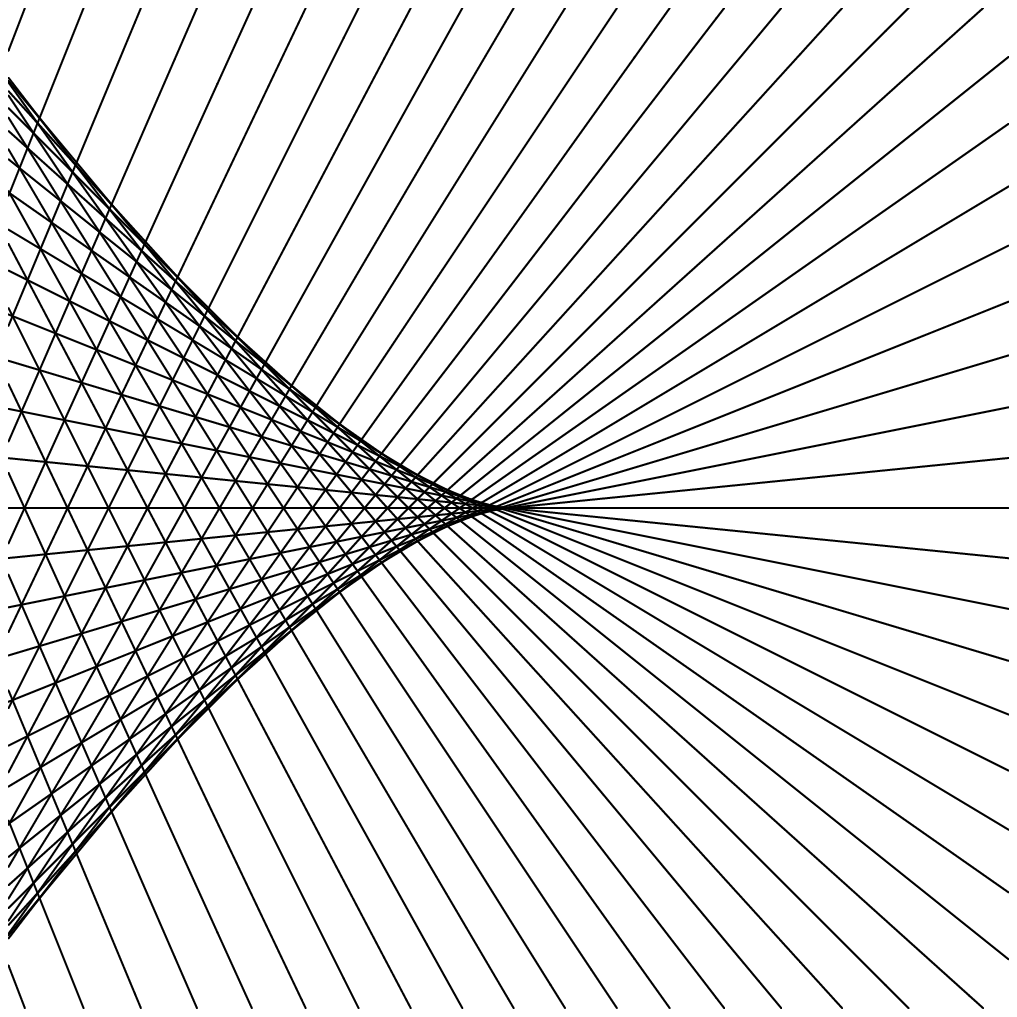} &
      \includegraphics[width=40mm]{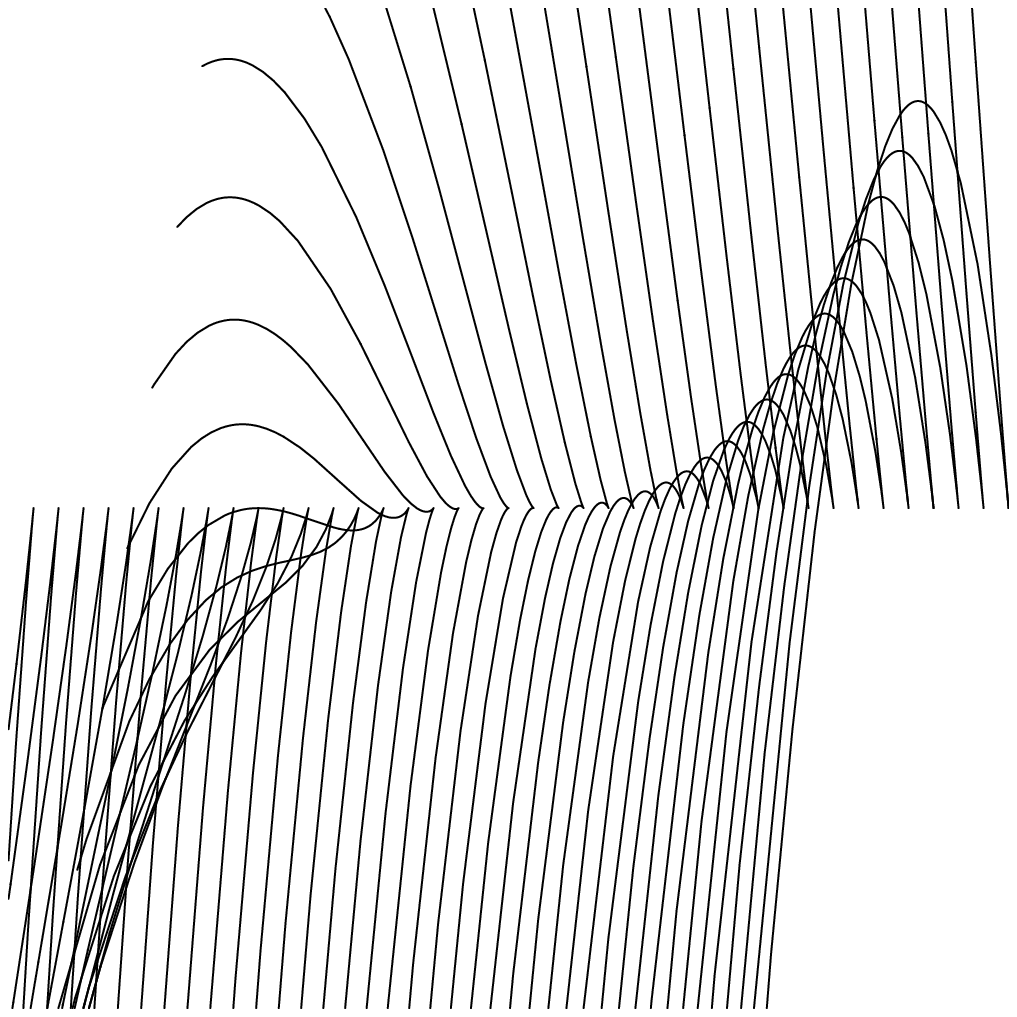} \\
Fig.5: (4) Regular cusp &
Fig.6:   (5)  Clairaut cusp &
 Fig.7:  (6)  Mixed fold \\

      \end{tabular}

  \end{center}

\end{figure}

\subsection{Quasi-linear first order partial differential equations}
We consider a time-dependent quasi-linear first order  partial differential equation
\[
\frac{\partial y}{\partial t}+\sum _{i=1}^n a_i(x,y,t)\frac{\partial y}{\partial x_i}
-b(x,y,t)=0,
\]
where $a_i(x,y,t)$ and $b(x,y,t)$ are $C^\infty$-function of $(x,y,t)=(x_1,\dots ,x_n,y,t).$
In order to clarify the situation in which there appeared a blow-up of the derivatives of solutions, we constructed a geometric framework
of the equation in \cite{Izu97}.
A {\it time-dependent quasi-linear first order partial differential equation} is defined by 
a hypersurface in $PT^*((\R^n\times \R)\times \R)$:
\[
E(1,a_1,\dots ,a_n,b)=\{(x,y,t),[ \xi :\eta:\sigma])\ |\ \sigma +\sum_{i=1}^n a_i(x,y,t)\xi_i
+b(x,y,t)\eta = 0 \ \}.
\]
A {\it geometric solution} of $E(1,a_1,\dots ,a_n,b)$ is a Legendrian submanifold $\mathscr{L}$ of $PT^*((\R^n\times \R)\times \R)$
lying in $E(1,a_1,\dots ,a_n,b)$ such that $\overline{\pi} |_{\mathscr{L}}$ is an embedding, where $\overline{\pi}:PT^*((\R^n\times \R)\times \R)\lon (\R^n\times \R)\times \R$ is the canonical projection.
Let $S$ be a smooth hypersurface in $(\R^n\times\R)\times \R$ . Then we have a unique Legendrian submanifold $\widehat{S}$ in 
$PT^*((\R^n\times \R)\times \R)$ such that $\overline{\pi}(\widehat{S})=S.$
It follows that if $\mathscr{L}$ is a geometric solution of $E(1,a_1,\dots ,a_n,b)$, then $\mathscr{L}=\widehat{\overline{\pi}(\mathscr{L})}.$
For any $(x_0,y_0,t_0)\in S,$ there exists a smooth submersion germ $f:((\R^n\times\R)\times \R,(x_0,y_0,t_0))\lon (\R,0)$
such that $(f^{-1}(0),(x_0,y_0,t_0))=(S,(x_0,y_0,t_0))$ as set germs.
A vector $\tau \partial/\partial t+\sum_{i=1}^n\mu _i\partial/\partial x_i+\lambda \partial/\partial y$ is tangent
to $S$ at $(x,y,t)\in (S,(x_0,y_0,t_0))$ if and only if $\tau \partial f/\partial t+\sum_{i=1}^n\mu _i\partial f /\partial x_i+\lambda \partial f/\partial y=0$ at $(x,y,t).$ Then we have the following representation of $\widehat{S}:$
\[
(\widehat{S},((x_0,y_0,t_0),[\sigma_0:\xi_0:\eta_0]))=\left\{\left((x,y,t),\left[\frac{\partial f}{\partial x}:\frac{\partial f}{\partial y}:\frac{\partial f}{\partial t}\right]\right) \Bigr | (x,y,t)\in (S,(x_0,y_0,t_0))\right\}.
\]
Under this representation, $\widehat{S}\subset E(1,a_1,\dots ,a_n,b)$ if and only if
\[
\frac{\partial f}{\partial t}+\sum_{i=1}^n a_i(x,y,t)\frac{\partial f }{\partial x_i}+b(x,y,t) \frac{\partial f}{\partial y}=0.
\]
Here, the {\it characteristic vector field} of $E(1,a_1,\dots ,a_n,b)$ is defined to be
\[
X(1,a_1,\dots ,a_n,b)=\frac{\partial }{\partial t}+\sum_{i=1}^n a_i(x,y,t)\frac{\partial  }{\partial x_i}+b(x,y,t) \frac{\partial }{\partial y}.
\]
In \cite{Izu97} a characterization theorem of geometric solutions was proved.
\begin{Th}[\cite{Izu97}]
Let $S$ be a smooth hypersurface in $(\R^n\times \R)\times \R$.
Then $\widehat{S}$ is a geometric solution of $E(1,a_1,\dots ,a_n,b)$ if and only if the characteristic vector field $X(1,a_1,\dots ,a_n,b)$ is tangent to $S.$
\end{Th}
\begin{Rem}{\rm We consider the Cauchy problem here:
\begin{eqnarray*}
&{}& \frac{\partial y}{\partial t}+\sum _{i=1}^n a_i(x,y,t)\frac{\partial y}{\partial x_i}
-b(x,y,t)=0, \\
&{}& y(0,x_1,\dots ,x_n)=\phi (x_1,\dots ,x_n),
\end{eqnarray*}
where $\phi$ is a $C^\infty$-function.
By Theorem 5.4, applying the classical method of characteristics, we can solve the above Cauchy problem.
Although $y$ is initially smooth, there is, in general, a critical time beyond which characteristics cross.
After the characteristics cross, the geometric solution becomes multi-valued.
Since the characteristic vector field $X(1,a_1,\dots ,a_n,b)$ is a vector field on the space $(\R^n\times \R)\times \R$,
the graph of the geometric solution $\overline{\pi}(\mathscr{L})\subset (\R^n\times \R)\times \R$ is a smooth hypersurface.
In general, however, $\widehat{\pi}_2|_{\overline{\pi}(\mathscr{L})}$ is a finite-to-one mapping, where $\widehat{\pi} _2:(\R^n\times \R)\times \R
\lon \R^n\times \R$ is $\widehat{\pi}_2 (x,y,t)=(x,t).$
}
\end{Rem}
The geometric solution $\mathscr{L}$ is a big Legendrian submanifold and it is Legendrian non-singular.
Therefore, the discriminant of the family of the momentary fronts
$\{W_t(\mathscr{L})\}_{t\in \R}$ is $\Delta_{\mathscr{L}} .$
We consider the following example:
\begin{eqnarray*}
&{}& \frac{\partial y}{\partial t}+2y\frac{\partial y}{\partial x}=0, \\
&{}& y(0,x)=\sin x,
\end{eqnarray*}
This equation is called {\it Burger's equation} and can be solved exactly by the method of characteristics.
We can draw the picture of the graph of the geometric solution  and the family of $\pi _2^{-1}(t)\cap W(\mathscr{L})$ in Fig.8.
\begin{figure}[htbp]
  \begin{center}
    \small
    \begin{tabular}{ccc}
       \includegraphics[width=50mm]{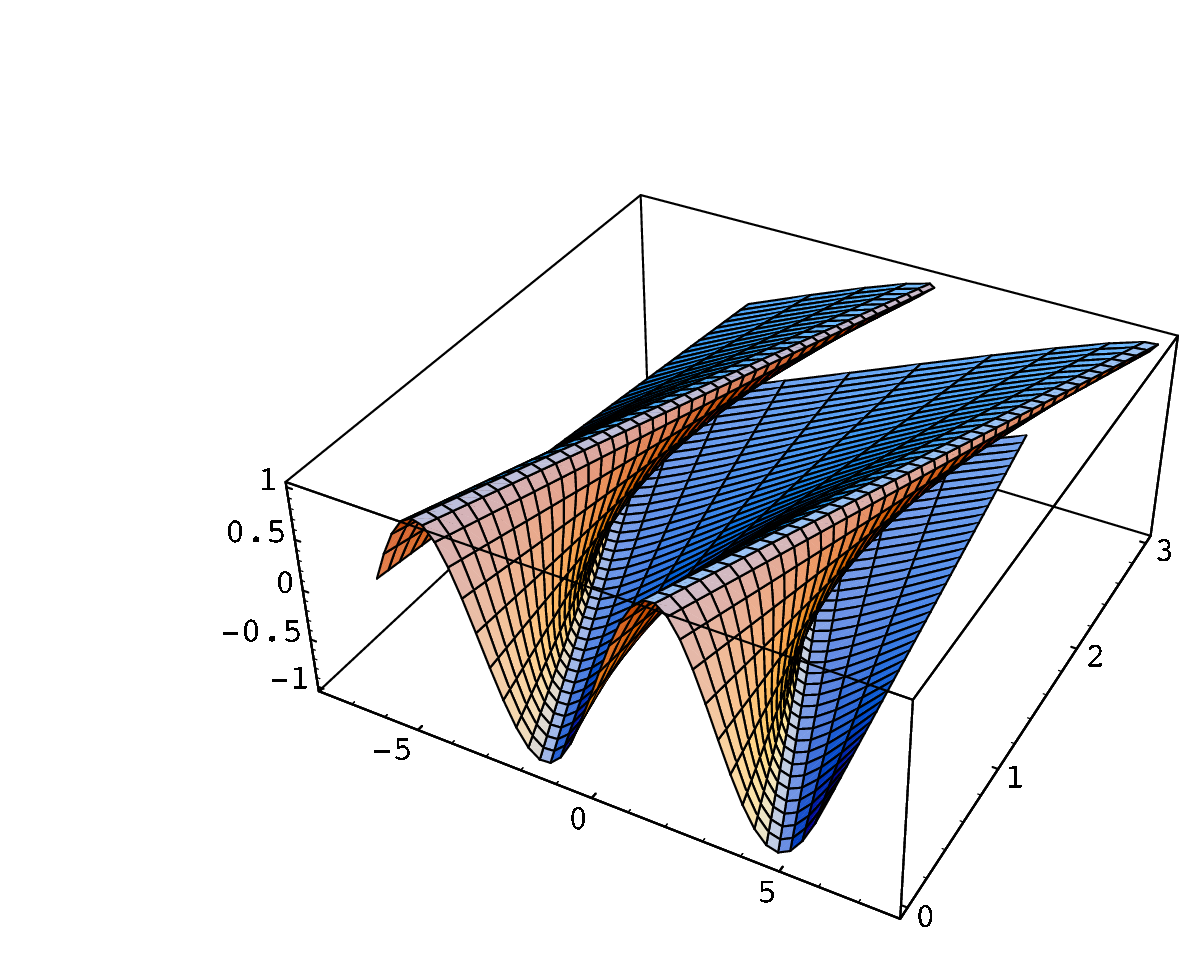} 
       & 
       {}
       &
       \includegraphics[width=40mm]{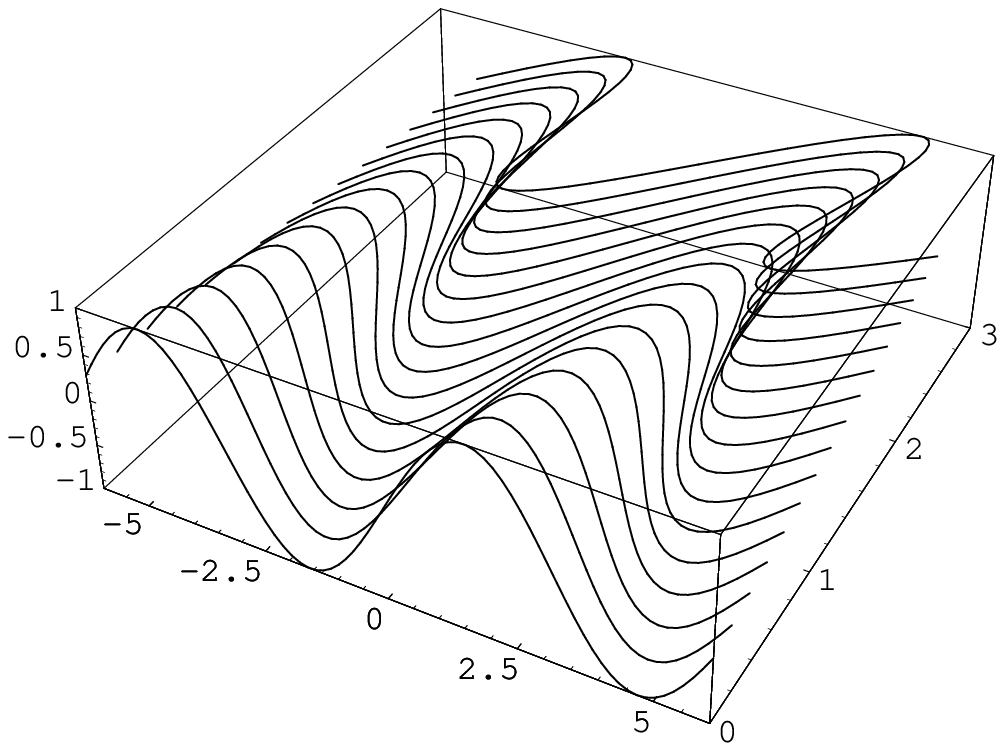}
     \\
The graph of the geometric solution \\
 of Burger's equation & {}& The family of $\pi _2^{-1}(t)\cap W(\mathscr{L})$ \\
 {} & Fig.8.& {}
 
 \\
      \end{tabular}

  \end{center}

\end{figure}
We can observe that the graph is a smooth surface in $(\R\times \R)\times \R$ but it is multi-valued.
Moreover, each $\pi _2^{-1}(t)\cap W(\mathscr{L})$ is non-singular but $\widehat{\pi}_2 |_{\pi _2^{-1}(t)\cap W(\mathscr{L})}$ has singularities.
Thus, $W(\mathscr{L})$ is a big wave front but not a graph-like wave front.

\subsection{Parallels and Caustics of hypersurfaces in Euclidean space}
In this subsection we respectively interpret the focal set (i.e., the evolute) of a hypersurface as the caustic and the parallels of a hypersurface as the graph-like 
momentary fronts by using the distance squared functions (cf. \cite{Izu07,Porteous}). 

Let $\bX:U \lon \R^n$ be an embedding, where $U$ is an open subset in $\R^{n-1}$. 
We denote that $M=\bX(U)$ and identify $M$ and $U$ via the embedding $\bX$. The {\it Gauss map} $\mathbb{G}:U \lon S^{n-1}$ is defined by $\mathbb{G}(u)=\bn(u)$, where $\bn(u)$ is the unit normal vector of $M$ at $\bX(u)$.
For a hypersurface $\bX:U \lon \R^n$, we define the {\it focal set} (or, {\it evolute}) of $\bX(U)=M$ by 
$$
{\rm F}_M=\bigcup _{i=1}^n\left\{\bX(u)+\frac{1}{\kappa_i(u)}\bn(u)\ |\ \kappa_i(u)\ {\rm is \ a\ principal \ curvature \ at \ }p=\bX(u), u \in U  \right\}
$$
and the set of {\it unfolded parallels} of $\bX(U)=M$ by 
$$
P_M=\left\{(\bX(u)+r \bn(u),r) \ |\ r \in \R \setminus \{0\}, u \in U \right\},
$$
respectively.
We also define the smooth mapping ${\rm F}_{\kappa_i}:U \lon \R^n$ and $P_{r}:U \lon \R^n$ by 
$$
{\rm F}_{\kappa_i}(u)=\bX(u)+\frac{1}{\kappa_i(u)}\bn(u), \ 
P_{r}(u)=\bX(u)+r \bn(u),
$$
where we fix a principal curvature $\kappa_i(u)$ on $U$ with $\kappa_i(u) \not=0$ and a real number $r \not=0$. 
\par
We now define families of functions in order to describe the focal set and the parallels of a hypersurface in $\R^n$. 
We define
$$
D:U \times (\R^n \setminus M) \lon \R
$$
by $D(u,\bv)=\| \bX(u)-\bv\|^2$ and 
$$
\overline{D}:U \times (\R^n \setminus M) \times \R_+ \lon \R
$$
by $\overline{D}(u,\bv,t)=\| \bX(u)-\bv\|^2-t$, where we denote that 
$\R_+$ is the set of positive real numbers. 
We call $D$ a {\it distance squared function} and $\overline{D}$ an {\it extended distance squared function} on $M=\bX(U)$. 
Denote that the function $d_{\bm{v}}$ and $\overline{d}_{\bm{v}}$ by $d_{\bm{v}}(u)=D(u,\bv)$ and $\overline{d}_{\bm{v}}(u,t)=\overline{D}(u,\bv,t)$ respectively.

By a straightforward calculation (cf., \cite{Izu07}), we have the following proposition:
\begin{Pro}
Let $\bX:U \lon \R^n$ be a hypersurface. Then 

$(1)$ $(\partial d_{\bm{v}}/\partial u_i)(u)=0 \ (i=1,\dots,n-1)$ if and only if there exists a real number $r \in \R \setminus \{0\}$ such that $\bv=\bX(u)+r\bn(u)$. 

$(2)$ $(\partial d_{\bm{v}}/\partial u_i)(u)=0 \ (i=1,\dots,n-1)$ and {\rm det}$(\mathcal{H}(d_v)(u))=0$ if and only if $\bv=\bX(u)+(1/\kappa(u)) \bn(u)$.

$(3)$ $\overline{d}_{\bm{v}}(u,t)=(\partial \overline{d}_{\bm{v}}/\partial u_i)(u,t)=0 \ (i=1,\dots,n-1)$ if and only if $\bv=\bX(u) \pm \sqrt{t}\bn(u).$

Here, $\mathcal{H}(d_{\bm{v}})(u)$ is the hessian matrix of the function $d_{\bm{v}}$ at $u$. 
\end{Pro}

As a consequence of Proposition 5.6, we have the following:
$$
C(D)=\left\{(u,\bv) \in U \times (\R^n \setminus M)\ |\ \bv=\bX(u)+r\bn(u),\ r \in \R \setminus \{0\}\right\},
$$
$$
\Sigma _*(\overline{D})=\left\{(u,\bv,t) \in U\times (\R^n \setminus M) \times \R_+\ |\ \bv=\bX(u) \pm \sqrt{t}\bn(u),\ u \in U \right\}.
$$
We can naturally interpret the focal set of a hypersurface as a caustic.
Moreover, the 
parallels of a hypersurface are given as a graph-like wave front (the momentary fronts). 
\begin{Pro}
For a hypersurface $\bX:U \lon \R^n$, the distance squared function 
$D:U \times (\R^n \setminus M) \lon \R$ is a Morse family of functions 
and the extended distance squared function $\overline{D}:U \times (\R^n \setminus M) \times \R_+ \lon \R$ is a non-degenerate graph-like Morse family of hypersurfaces.
\end{Pro}
\demo
By Proposition 4.5, it is enough to show that $D$ is a Morse family of hypersurfaces.
For any $\bv=(v_1,\dots,v_n) \in \R^n \setminus M$, we have 
$D(u,\bv)=\sum^{n}_{i=1}(x_i(u)-v_i)^2$, where $\bX(u)=(x_1(u),\dots,x_n(u))$. 
We shall prove that the mapping 
$$
\Delta ^*D=\left(D, \frac{\partial D}{\partial u_1},\dots,\frac{\partial D}{\partial u_{n-1}} \right)
$$
is a non-singular at any point. 
The Jacobian matrix of $\Delta ^*D$ is given by 
$$\left(
\begin{array}{cccccc}
A_1(u)  &\cdots &A_{n-1}(u) &-2(x_1(u)-v_1)&\cdots &-2(x_{n-1}-v_{n-1})  \\
A_{11}(u)&\cdots&A_{1(n-1)}(u)&-2x_{1u_1}(u)&\cdots&-2x_{nu_1}(u)\\
\vdots&\ddots&\vdots&\vdots&\ddots&\vdots\\
A_{(n-1)1}(u)&\cdots&A_{(n-1)(n-1)}&-2x_{1u_{n-1}}(u)&\cdots&-2x_{nu_{n-1}}(u)
\end{array}
\right),
$$
where $A_i(u)=\langle 2\bX _{u_i}(u),\bX(u)-\bv\rangle$, $A_{ij}(u)=2(\langle \bX_{u_iu_j}(u), \bX(u)-\bv \rangle+\langle \bX_{u_i}(u), \bX_{u_j}(u) \rangle)$ and $\langle , \rangle$ is the inner product of $\R^n$. 
Suppose that $(u,\bv,t_0)\in \Sigma _*(\overline{D}).$ Then we have $\bv=\bX(u)\pm\sqrt{t_0}\bn(u).$
Therefore, we have
\[
J_{\Delta ^*D}(u,\bv,t_0)=\left(
\begin{array}{ccc}
0&\mp 2\sqrt{t_0}\bn(u) \\
A_{ij}(u)&-2\bX_{u_i}(u)
\end{array}
\right).
\]
Since $\bn(u),\bX_{u_1}(u),\dots ,\bX_{u_{n-1}}(u)$ are linearly independent,
the rank of $J_{\Delta ^*D}(u,\bv,t_0)$ is $n.$
This means that $D$ is a Morse family of hypersurfaces. 
\enD

By the method for constructing a Lagrangian submanifold germ from a Morse family of functions (cf. \S 2), we can define a Lagrangian submanifold germ whose generating family is the distance squared function $D$ of $M=\bX(U)$ as follows: 
For a hypersurface $\bX:U \lon \R^n$ where $\bX(u)=(x_1(u),\dots,x_n(u))$, we define 
$$
L(D):C(D) \lon T^* \R^n
$$
by
$$
L(D)(u,\bv)=(\bv,-2(x_1(u)-v_1),\dots,-2(x_n(u)-v_n)),
$$
where $\bv=(v_1,\dots,v_n)$.
\par
On the other hand, by the method for constructing the graph-like Legendrian unfolding from a graph-like Morse family of hypersurfaces (cf. \S 4), we can define a graph-like Legendrian unfolding whose generating family is the extended distance squared function $\overline{D}$ of $M=\bX(U)$.  
For a hypersurface $\bX:U \lon \R^n$ where $\bX(u)=(x_1(u),\dots,x_n(u))$, we define 
$$
{\mathfrak{L}}_{D}:C(D) \lon J^1_{GA}(\R^n,\R)
$$
by 
$$
{\mathfrak{L}}_{D}(u,\bv)=(\bv,\| \bX(u)-\bv \|^2,-2(x_1(u)-v_1),\dots,-2(x_n(u)-v_n)),
$$
where $\bv=(v_1,\dots,v_n)$.
\begin{Co}Using the above notation, $L(D)(C(D))$ is a Lagrangian submanifold such that the distance squared function $D$ is the generating family of $L(D)(C(D))$ and 
${\mathfrak{L}}_{D}(C(D))$ is a non-degenerate graph-like Legendrian unfolding such that 
the extended distance squared function $\overline{D}$ is the graph-like generating family of ${\mathfrak{L}}_{D}(C(D))$.
\end{Co}
\par
By Proposition 5.6, the caustic $C_{L(D)(C(D))}$ of $L(D)(C(D))$ is the focal set $F_M$
and the graph-like wave front $W({\mathfrak{L}}_{D}(C(D)))$ is the set of unfolded parallels $P_M.$
\par
We now briefly describe the theory of contact with foliations. 
Here we consider the relationship between the contact of submanifolds with foliations and 
the ${\mathcal R}^+$-class of functions.  
Let $X_i\ (i=1,2)$ be submanifolds of $\R^n$ with ${\rm dim}\, X_1={\rm dim}\, X_2$, 
$g_i:(X_i,\bar{x}_i) \lon (\R^n,\bar{y}_i)$ be immersion germs and 
$f_i:(\R^n,\bar{y}_i) \lon (\R,0)$ be submersion germs. 
For a submersion germ $f: (\R^n,0) \lon (\R,0)$, we have  
 the regular foliation $\mathcal{F}_{f}$ defined by $f$; i.e.,
$\mathcal{F}_{f}=\{ f^{-1} (c) | c \in (\R,0) \}$.
We say that {\it the contact of $X_1$ with the regular foliation $\mathcal{F}_{f_1}$ } 
at $\bar{y}_1$ is of the {\it same type} as {\it the contact of 
$X_2$ with the regular foliation $\mathcal{F}_{f_2}$} at $\bar{y}_2$
if there is a diffeomorphism germ
$\Phi:(\R^n,\bar{y}_1) \lon (\R^n,\bar{y}_2)$
 such that $\Phi(X_1)=X_2$ and
$\Phi(Y_1 (c)) = Y_2 (c)$, where $Y_i (c) = f^{-1}_i (c)$
for each $c \in (\R,0)$. 
In this case we write 
$K(X_1,\mathcal{F}_{f_1};\bar{y}_1) = K(X_2,\mathcal{F}_{f_2};\bar{y}_2)$.
 We apply the method of Goryunov\cite{Go} to 
 the case for ${\mathcal R}^+$-equivalences among function germs. Then we 
 have the following proposition: 
\begin{Pro} {\rm (\cite[Appendix]{Go})}
Let $X_i\ (i=1,2)$ be submanifolds of $\R^n$ with ${\rm dim}\, X_1 = {\rm dim}\, X_2=n-1$
{\rm (}i.e. hypersurfaces{\rm )}, 
$g_i:(X_i,\bar{x}_i) \lon (\R^n,\bar{y}_i)$ be immersion germs and 
$f_i:(\R^n,\bar{y}_i) \lon (\R,0)$ be submersion germs.
Then $K(X_1,\mathcal{F}_{f_1};\bar{y}_1) = K(X_2,\mathcal{F}_{f_2};\bar{y}_2)$ 
if and only if 
$f_1 \circ g_1$ and $f_2 \circ g_2$ are ${\mathcal R}^+$-equivalent.
\end{Pro}
\par
\par
On the other hand, we define a function ${\mathcal D}:\R^n\times\R^n\lon \R$ by
${\mathcal D}(\bx,\bv)=\|\bx-\bv\|^2.$
For any $\bv \in \R^n\setminus M$, we write 
$\mathfrak{d}_{\bm{v}}(\bx) = {\mathcal D} (\bx,\bv)$ and 
we have a hypersphere ${\mathfrak{d}_{\bm{v}}}^{-1}(c ) = S^{n-1} (\bv,\sqrt{c})$
for any $c>0.$ 
It is easy to show that $\mathfrak{d}_{\bm{v}}$ is a submersion.
For any $u \in U$, we consider 
$\bv^\pm = \bX(u) \pm \sqrt{c}\bn(u)\in \R^n\setminus M$. Then we have 
$$
\mathfrak{d}_{\bm{v}^\pm} \circ \bX (u) = \mathcal{D} \circ 
(\bX \times id_{\R^n})(u,\bv^\pm)=c,
$$
and
$$
\frac{\partial{\mathfrak{d}_{\bm{v}^\pm}} \circ \bX}{\partial{u_i}}(u)
=\frac{\partial{D}}{\partial{u_i}}(u,\bv^\pm)=0.
$$
for $i=1,\dots,n-1$. 
This means that the hyperspheres $\mathfrak{d}_{\bm{v}^\pm}^{-1} (r) =
S^{n-1}(\bv^\pm,\sqrt{c})$ are tangent to $M=\bX(U)$ at $p=\bX(u)$. 
In this case, we call each one of 
 $S^{n-1}(\bv^\pm,\sqrt{c})$ a {\it tangent hypersphere} at $p=\bX(u)$
with the center $\bv ^\pm.$ 
However, there are infinitely many tangent hyperspheres at a general point
$p=\bX (u)$ depending on the real number $c .$
If $\bv$ is a point of the focal set (i.e., $\bv=F_{\kappa}(u)$ for some $\kappa$), the tangent hypersphere with the center 
$\bv$ is called the {\it osculating hypersphere} (or, {\it curvature hypersphere}) at $p=\bX (u)$
which is uniquely determined. 
For $\bv^\pm = \bX(u) \pm \sqrt{c}\bn(u),$ we also have  regular foliations
\[
{\mathcal F}_{\mathfrak{d}_{\bm{v}^\pm}}=\left\{S^{n-1}(\bv^\pm ,\sqrt{t})\ \Bigm |\ t\in \left(\R,c\right)\right\}
\]
whose leaves are hyperspheres with the center $\bv^\pm$ such that the case $t=c$ 
corresponds to the tangent hypersphere with radius $|c|.$
Moreover, if $\bv=F_\kappa (u),$ then $S^{n-1}(\bv ,1/\kappa (u))$
is the osculating hypersphere. In this case $(\bX^{-1}({\mathcal F}_{\mathfrak{d}_{\bm{v}}}),u)$ is a singular foliation germ at $u$
which is called an {\it osculating hyperspherical foliation} of $M=\bX(U)$ at $p=\bX(u)$ (or, $u$).
We denote it by ${\mathcal O\mathcal F}(M,u).$
Moreover, if $\bv\in M_{{\mathfrak{L}}_{D}(C(D))}$, then there exists $r_0\in \R\setminus \{0\}$ such that $(\bv,r_0)$ is a self-intersection point of $P_M,$
so that there exist different $u,v \in U$
such that $\bv=\bX(u)+r_0\bn (u)=\bX(v)+r_0\bn(v).$
Therefore, the hypersphere $S^{n-1}(\bv,|r_0|)$ is tangent to $M=\bX(U)$ at both the points $p=\bX(u)$ and $q=\bX(v).$
Then we have an interpretation of the geometric meanings of the Maxwell stratified set in this case:
\[
M_{{\mathfrak{L}}_{D}(C(D))}=\{\bv \ |\ \exists r_0\in \R\setminus \{0\}, S^{n-1}(\bv,|r_0|)\ \mbox{is\ tangent\ to}\ M\ \mbox{at\ least\ two\ different\ points}\}.
\]
Therefore, we call the Maxwell stratified set $M_{{\mathfrak{L}}_{D}(C(D))}$ the {\it set of the centers of multiple tangent spheres} of $M.$
\par
We consider the contact of hypersurfaces with families of hyperspheres.
Let $\bX_i : (U,\bar{u}_i) \lon (\R^n,p_i)$ $(i=1,2)$ be hypersurface germs.
We consider distance squared functions 
$D_i :(U \times \R^n,(\bar{u}_i,\bv_i)) \lon \R$ of $M_i=\bX_i(U)$, 
where $\bv_i={\rm Ev}_{\kappa _i}(\bar{u}_i).$ 
We write $d_{i,\bm{v}_i} (u) = D_i (u,\bv_i)$, then we have 
$d_{i,\bm{v}_i}(u)=\mathfrak{d}_{\bm{v}_i} \circ \bX_i (u)$. 
Then we have the following theorem:
\begin{Th}
Let $\bX_i : (U,\bar{u}_i) \lon \R^n,p_i)$ $(i=1,2)$ be hypersurface germs
 such that the corresponding graph-like  Legendrian unfolding germs
 $\mathfrak{L}_{D_i}(C(D_i))$ are $S.P^+$-Legendrian stable (i.e., the corresponding Lagrangian
 submanifold germs  
$L(D_i) (C(D_i))$ are 
Lagrangian stable), 
where $\bv_i={\rm Ev}_{\kappa _i}(\bar{u}_i)$ are centers of the osculating hyperspheres of $M_i=\bX _i(U)$ respectively. 
Then the following conditions are equivalent:
\par\noindent
{\rm (1)} $\mathfrak{L}_{D_1}(C(D_i))$ and $\mathfrak{L}_{D_i}(C(D_2))$ are $S.P^+$-Legendrian equivalent,
\par\noindent
{\rm (2)} $\overline{D}_1$ and $\overline{D_2}$ are $s$-$S.P^+$-$\mathcal{K}$-equivalent,
\par\noindent
{\rm (3)} $\overline{d}_{1,\bm{v}_1}$ and $\overline{d}_{1,\bm{v}_2}$ are $S.P$-$\mathcal{K}$-equivalent,
\par\noindent
{\rm (4)} $d_{1,\bm{v}_1}$ and $d_{2,\bm{v}_2}$ are $\mathcal{R}$-equivalent, 
\par\noindent
{\rm (5)}  $K(M_1, \mathcal{F}_{\mathfrak{d}_{\bm{v}_1}};p_1)=
K(M_2 , \mathcal{F}_{\mathfrak{d}_{\bm{v}_2}};p_2)$,
\par\noindent
{\rm (6)} $D_1$ and $D_2$ are $P$-$\mathcal{R^+}$-equivalent,
 \par\noindent
{\rm (7)} $L(D_1)(C(D_1))$ and $L(D_2)(C(D_2))$ are Lagrangian equivalent,
\par\noindent
{\rm (8)} $P_{M_1}$ and $P_{M_2}$ are $S.P^+$-diffeomorphic.
\end{Th}
\demo
By Theorem 5.9, the conditions (4) and (5) are equivalent. 
By the assertion (3) of Proposition 5.6, we have
$W({\mathfrak{L}}_{D_i}(C(D_i))=P_{M_i}.$
Thus, the other conditions are equivalent to each other by Theorem 4.4.
\enD

\par
We remark that if $L(D_1)$ and $L(D_2)$ are Lagrangian equivalent, then the
corresponding caustics are diffeomorphic.
Since the caustic of $L(D)$ is the focal set of a hypersurface $M=\bX (U),$ the above theorem gives a symplectic interpretation for the contact of hypersurfaces with
family of hyperspheres.
Moreover, the $S.P^+$-diffeomorphism between the graph-like wave front sets sends the Maxwell stratified sets to each other.
Therefore,
we have the following corollary.
\begin{Co} Under the same assumptions as those of the above theorem
for hypersurface germs $\bX_i : (U,\bar{u}_i) \lon (\R^n,p_i)$ $(i=1,2)$, we have
the following:
If one of the conditions of the above theorem is satisfied, then
\par
{\rm (1)} The focal sets $F_{M_1}$ and $F_{M_2}$
are diffeomorphic as set germs.
\par
{\rm (2)} The osculating hyperspherical foliation germs
${\mathcal O\mathcal F}(M_1,\bar{u}_1)$, ${\mathcal O\mathcal F}(M_2,\bar{u}_2)$
are diffeomorphic.
\par
{\rm (3)} The sets of the centers of multiple tangent spheres 
of $M_1$ and $M_2$
are diffeomorphic as set germs.
\end{Co}

\subsection{Caustics of world sheets}
 Recently the author has discovered an application of the theory of graph-like Legendrian unfoldings to 
the caustics of world sheets in Lorentz space forms. 
In the theory of relativity, we do not have the notion of 
time constant, so that everything that is moving depends on the time.
Therefore, we have to consider world sheets instead of spacelike submanifolds.
Let $\mathbb{L}^{n+1}_1$ be an $n+1$-dimensional Lorentz space form (i.e., Lorentz-Minkowski space, de Sitter space or anti-de Sitter space). For basic
concepts and properties of Lorentz space forms, see \cite{Oneil}.
We say that a non-zero vector $\bx \in \mathbb{L}^{n+1}_1$ is
{\it spacelike, lightlike or timelike} if $\langle \bx,\bx \rangle
>0,$ $\langle \bx,\bx \rangle =0$ or  $\langle \bx,\bx \rangle <0$,
respectively. 
Here, $\langle \bx,\by \rangle$  is the induced pseudo-scalar product of $\mathbb{L}^{n+1}_1.$
We only consider the local situation here.
Let $\bX :U\times I\lon \mathbb{L}^{n+1}_1$ be a timelike embedding of codimension $k-1,$
where $U\subset \R^s$ ($s+k=n+1$) is an open subset and $I$ an open interval.
 We write 
$W=\bX(U\times I)
$ and identify $W$ and $U\times I$ via the embedding $\bX.$
Here, the embedding $\bX$ is said to be {\it timelike} if the tangent space $T_p  W$
of $W$ at $p=\bX(u,t)$ is a timelike subspace (i.e., Lorentz subspace of $T_p\mathbb{L}_1^{n+1}$) for any point $p\in W$.
We write $\mathcal{S}_t=\bX(U\times\{t\})$ for each $t\in I.$
We call $\mathscr{F}_S=\{\mathcal{S}_t\ |t\in I\}$ a {\it spacelike foliation} 
on $W$ if $\mathcal{S}_t$ is a spacelike submanifold for any $t\in I.$
Here, we say that $\mathcal{S}_t$ is {\it spacelike} if the tangent space $T_p\mathcal{S}_t$ 
consists only spacelike vectors (i.e., spacelike subspace) for
any point $p\in \mathcal{S}_t.$
We call $\mathcal{S}_t$ a {\it momentary space} of $\mathscr{F}_S=\{\mathcal{S}_t\ |t\in I\}$.
We say that $W=\bX(U\times I)$ (or, $\bX$ itself) is a {\it world sheet}
if $W$ is time-orientable.
It follows that there exists a unique timelike future directed unit normal vector field $\bn^T(u,t)$ along $\mathcal{S}_t$ on $W$ (cf., \cite{Oneil}). This means that $\bn^T(u,t)\in T_pW$ and is pseudo-orthogonal to $T_p\mathcal{S}_t$ for
$p=\bX(u,t).$
Since $T_pW$ is a timelike subspace of $T_p\mathbb{L}^{n+1}_1,$
the pseudo-normal space $N_p(W)$ of $W$ is a $k-1$-dimensional spacelike subspace of $T_p\mathbb{L}^{n+1}_1$
(cf.,\cite{Oneil}). On the pseudo-normal space $N_p(W),$ we have a $(k-2)$-unit sphere
\[
N_1(W)_p=\{\bxi\in N_p(W)\ |\ \langle \bxi,\bxi\rangle =1\ \}.
\]
Therefore, we have a unit spherical normal bundle over $W$:
\[
N_1(W)=\bigcup _{p\in W} N_1(W)_p.
\]
For an each momentary space $\mathcal{S}_t,$ we have a unit spherical normal bundle $N_1[\mathcal{S}_t]=N_1(W)|_{\mathcal{S}_t}$
over $\mathcal{S}_t.$
Then we define a hypersurface
$
\mathbb{LH}_{\mathcal{S}_t}:N_1[\mathcal{S}_t]\times \R\lon \mathbb{L}_1^{n+1}
$
by
$$
\mathbb{LH}_{\mathcal{S}_t}(((u,t),\bxi),\mu)=\bX(u,t)+\mu (\bn^T(u,t)+\bxi),
$$
where $p=\bX (u,t)$, which is called the {\it lightlike hypersurface\/} in the Lorentz space form $\mathbb{L}_1^{n+1}$ along $\mathcal{S}_t$.
The lightlike hypersurface of a spacelike submanifold in a Lorentz space form has been defined and investigated in
\cite{Izu-Sato1,Izu-Sato2,Izu-Sato3}.
The set of singular values of the lightlike hypersurface is called a {\it focal set} of $\mathcal{S}_t$.
We remark that the situation is different from the Riemannian case. The lightlike hypersurface is a wave front in $\mathbb{L}^{n+1}_1$,
so that the focal set is the set of Legendrian singular values. In the Riemannian case, the focal set is the set of Lagrangian singular values. In the Lorentzian case, we consider world sheets instead of a single spacelike submanifold.
Since a world sheet is a one-parameter family of spacelike submanifolds, we can naturally apply the theory of wave front propagations.
We define
\[
\widetilde{\mathbb{LH}}(W)=\bigcup _{t\in I} \mathbb{LH}_{\mathcal{S}_t}(N_1[\mathcal{S}_t]\times \R)\times \{t\}
\subset \mathbb{L}^{n+1}_1\times I,
\]
which is called a {\it unfolded lightlike hypersurface}.
In \cite{Izu14-2} we show that the unfolded lightlike hypersurface is a graph-like wave front and each lightlike hypersurface
is a momentary front for the case that $\mathbb{L}^{n+1}_1$ is the anti-de Sitter space.
One of the motivations for investigating this case is given in the brane world scenario (cf., \cite{Bousso-Randal,Bousso}).
There, lightlike hypersurfaces and caustics along world sheets have been considered in the simplest case.
Since the unfolded lightlike hypersurface is a graph-like Legendrian unfolding, we can investigate
not only the caustic but also the Maxwell stratified set as an application of the theory of Legendrian unfoldings.
We can apply Theorem 4.1 to this case and get some geometric information on world sheets.
We can also consider the {\it lightcone pedal} of world sheets and investigate the geometric properties as an application of the theory of graph-like unfoldings \cite{IzuSM,Izu14}.

\subsection{Control theory}
In \cite{Zakalyukin95,Zakalyukin05} Zakalyukin applied $S.P^+$-Legendrian equivalence to the study of 
 problems which occur in the control theory.  In \cite{Zakalyukin95} he has given the following simple example: Consider a plane $\R^2$. For each point $q=(q_1,q_2)\in \R^2,$ we consider an admissible curve on the tangent plane $\R^2=T_q\R^2$ defined by $p_1=1+u, p_2=u^2\ (u\in\R)$,
 where $(p_1,p_2)\in \R^2$ is the coordinates of $\R^2=T_q\R^2$.
 So this admissible curve is independent of the base point $q\in \R^2.$
 The initial front is given by $W_0=\{(q_1,f(q_1))\ |\ q_1\in \R\}$ for some function $f(q_1).$
 According to the Pontryagin maximum principle, externals of the corresponding time optimal control problem are defined by a
 canonical system of equations with the Hamiltonian $H(p,q)=\max _u (p_1(1+u)+p_2u^2).$
 This system can be solved exactly and the corresponding family of fronts $W_t$ are given parametrically in the form $W_t=\Phi _t(W_0)$:
 \[
 \Phi _t(q_1,t)=\left(q_1+t\left(1+\frac{1}{2}\frac{df}{dq_1}\right),f(q_1)+\frac{t}{4}\left(\frac{df}{dq_1}\right)^2\right).
 \]
 Under the condition $f'(0)=0$ and $f''(0)>0$, he has shown that the picture of the  discriminant set of the family $\{W_t\}_{t\in I}$ is the same as that of the discriminant set of the germ (6) of Theorem 5.2.
 He also applies $S.P^+$-Legendrian equivalence to translation-invariant control problems in \cite{Zakalyukin05}.
 \par
 The author is not a control theory specialist, so that he cannot explain the results in detail here.
 However, it seems that there might be a lot of applications of the theory of wave front propagations to this area.
 For the detailed arguments, see the original articles.

\begin{flushleft}
\textsc{
Shyuichi Izumiya
\\ Department of Mathematics
\\ Hokkaido University
\\ Sapporo 060-0810, Japan} 
\par
e-mail: izumiya@math.sci.hokudai.ac.jp
\end{flushleft}

\begin{thebibliography}{99999}
{\renewcommand{\baselinestretch}{1}
\small
\bibitem{Arnold1}
V. I. Arnol'd, S. M. Gusein-Zade and A. N. Varchenko,
\newblock{\em Singularities of Differentiable Maps vol. I}.
\newblock Birkh\"auser, 1986.

\bibitem{Arnold-blue}
V. I. Arnol'd,
\newblock{\em Contact geometry and wave propagation}.
\newblock Monograph. Enseignement Math. 34 (1989)

\bibitem{Arnold-pink}
V. I. Arnol'd,
\newblock{\em Singularities of caustics and wave fronts}.
\newblock Math. Appl. 62, Kluwer , Dordrecht, 1990.

\bibitem{Bousso-Randal}
R. Bousso  and L Randall, 
\newblock{\em Holographic domains of ant-de Sitter space}, 
\newblock Journal of High Energy Physics,04 (2002), 057
\bibitem{Bousso}
R. Bousso,
\newblock{\em The holographic principle},
\newblock REVIEWS OF MODERN PHYSICS 74 (2002), 825--874.
\bibitem{Bro}
TH. Br\"ocker,
\newblock{\em Differentiable Germs and Catastrophes}.
\newblock London Mathematical Society Lecture Note Series 17, Cambridge University Press, 1975.

\bibitem{Bruce}
J. W. Bruce,
\newblock{\em Wavefronts and parallels in Euclidean space}.
\newblock Math. Proc. Cambridge Philos. Soc. {\bf 93} (1983) 323--333
\bibitem{Damon}
J. Damon,
\newblock{The unfolding and determinacy theorems for subgroups of $\mathcal{A}$ and
$\mathcal{K}$}.
\newblock Memoirs of A.M.S. {\bf 50} No. {\bf 306}, (1984)


\bibitem{Dimca}
A. Dimca,
\newblock{\em Topics on real and complex singularities}. 
\newblock Advanced Lectures in Mathematics (1987)

\bibitem{Go}
V. V. Goryunov,
{\em Projections of Generic Surfaces with Boundaries},
Adv. Soviet Math., {\bf 1} (1990), 157--200


\bibitem{Gor-Zak}
V. Goryunov and V. M. Zakalyukin,
\newblock{Lagrangian and Legendrian Singularities}.
\newblock Real and Complex Singularities, Trends in Mathematics, 169--185, Birkh\"auger, 2006 

\bibitem{HIIY}
A. Hayakawa, G. Ishikawa, S.Izumiya and K. Yamaguchi,
\newblock
{\em Classification of generic integral diagrams and first order ordinary differential equations}.
\newblock International Journal of Mathematics,  {\bf 5} (1994), 447--489.

\bibitem{Hormander}
L. H\"ormander,
\newblock{\em Fourier Integral Operators,I}.
\newblock Acta. Math. {\bf 128} (1972), 79--183
\bibitem{Izudoc}
S. Izumiya, 
\newblock{Generic bifurcations of varieties}.
\newblock manuscripta math. {\bf 46} (1984), 137--164

\bibitem{Izumiya93}
S. Izumiya,
\newblock{\em Perestroikas of optical wave fronts and graphlike Legendrian unfoldings}.
\newblock J. Differential Geom. {\bf 38} (1993), 485--500.

\bibitem{Izu95}
S. Izumiya,
\newblock{\em Completely integrable holonomic systems of first-order differential equations}.
\newblock Proc. Royal Soc. Edinburgh 125A (1995), 567--586.

\bibitem{Izu-Kro}
S. Izumiya and Y. Kurokawa,
\newblock{\em Holonomic systems of Clairaut type}.
\newblock Differential Geometry and its Applications 5 (1995), 219--235.


\bibitem{Izu97}
S. Izumiya and G. Kossioris,
\newblock{\em Geometric Singularities for Solutions of Single Conservation Laws}.
\newblock Arch. Rational Mech. Anal.  139 (1997), 255--290.



\bibitem{Izu07}
S. Izumiya,
\newblock{\em Differential Geometry from the viewpoint of Lagrangian or Legendrian singularity theory}.
\newblock in {\it Singularity Theory, Proceedings of the 2005 Marseille Singularity School and Conference, by D. Ch\'eniot et al}. World Scientific (2007) 241--275.

\bibitem{IzuSM} S. Izumiya and M.C. Romero Fuster, {\em The lightlike flat geometry
on spacelike submanifolds of codimension two in Minkowski space.} Selecta Math. (N.S.) 13 (2007), no. 1, 23--55.

\bibitem{Izu-Sato1}
S. Izumiya and T. Sato,
\newblock{\em Lightlike hypersurfaces along spacelike submanifolds in Minkowski space-time}.
\newblock{ Journal of Geometry and Physics}. {\bf 71} (2013), 30--52.

\bibitem{Izu-Sato2}
S. Izumiya and T. Sato,
\newblock{\em Lightlike hypersurfaces along spacelike submanifolds in anti-de Sitter space},
preprint (2012)
\bibitem{Izu-Sato3}
S. Izumiya and T. Sato,
\newblock{\em Lightlike hypersurfaces along spacelike submanifolds in de Sitter space},
preprint (2013)


\bibitem{Izumiya-Takahashi} S. Izumiya and M. Takahashi,
\newblock{\em Spacelike parallels and evolutes in Minkowski pseudo-spheres}.
\newblock{Journal of Geometry and Physics}. {\bf 57} (2007), 1569--1600.

\bibitem{Izumiya-Takahashi2} S. Izumiya and M. Takahashi,
\newblock{\em Caustics and wave front propagations: Applications to differential geometry}.
\newblock{Banach Center Publications. Geometry and topology of caustics}. {\bf 82} (2008)  125--142.

\bibitem{Izumiya-Takahashi3} S. Izumiya and M. Takahashi,
\newblock{\em Pedal foliations and Gauss maps of hypersurfaces in Euclidean space}.
\newblock{Journal of Singularities}. {\bf 6} (2012)  84--97.

\bibitem{Izu14}
S. Izumiya,
\newblock{\em Geometry of world sheets 
in
Lorentz-Minkowski space},
preprint (2014).


\bibitem{Izu14-2}
S. Izumiya,
\newblock{\em Caustics of world sheets in anti-de Sitter space},
in preparation.


\bibitem{Martinet}
J. Martinet,
\newblock{Singularities of Smooth Functions and Maps}.
\newblock London Math. Soc. Lecture Note Series, Cambridge University Press {\bf 58} (1982)

\bibitem{Maslov}
V. P. Malsov and M. V. Fedoruik,
\newblock{\em Semi-classical approximation in quantum mechanics},
\newblock D. Reidel Publishing Company, Dordrecht, 1981.

\bibitem{Oneil}
B. O'Neill,
\newblock{\em Semi-Riemannian Geometry},
\newblock Academic Press, New York, 1983.

\bibitem{Porteous}
I. Porteous,
\newblock{\em The normal singularities of submanifold}.
\newblock J. Diff. Geom. {5}, (1971), 543--564.

\bibitem{Zak}
V. M. Zakalyukin,
\newblock{\em Lagrangian and Legendrian singularities},
\newblock Funct. Anal. Appl. (1976), 23--31.

\bibitem{Zak76}
V. M. Zakalyukin,
\newblock{\em Reconstructions of fronts and caustics depending one
parameter},
\newblock Funct. Anal. Appl. (1976), 139--140.


\bibitem{Zak1}
V. M. Zakalyukin,
\newblock{\em Reconstructions of fronts and caustics depending one
parameter and versality of mappings},
\newblock J. Sov. Math. 27 (1984), 2713--2735.


\bibitem{Zakalyukin95} 
V.M. Zakalyukin, 
\newblock{\em Envelope of Families of Wave Fronts and Control Theory.}
\newblock Proc. Steklov Inst. Math. {\bf 209} (1995), 114--123.

\bibitem{Zakalyukin05} 
V.M. Zakalyukin, 
\newblock{\em  Singularities of Caustics in generic translation-invariant control problems.}
\newblock Journal of Mathematical Sciences, {\bf 126} (2005), 1354--1360.


}
\end{thebibliography}
\end{document}